# Minimal and minimal invariant Markov bases of decomposable models for contingency tables

HISAYUKI HARA[1], SATOSHI AOKI[2] and AKIMICHI TAKEMURA[3]

[1]*Department of Technology Management for Innovation, University of Tokyo, 7-3-1 Hongo, Bunkyo-ku, Tokyo 113-8656, Japan. E-mail: hara@tmi.t.u-tokyo.ac.jp*
[2]*Graduate School of Science and Engineering (Science Course), Kagoshima University 1-21-35 Korimoto, Kagoshima-shi, Kagoshima 890-0065, Japan. E-mail: aoki@sci.kagoshima-u.ac.jp*
[3]*Graduate School of Information Science and Technology, University of Tokyo, 7-3-1 Hongo, Bunkyo-ku, Tokyo 113-8656, Japan. E-mail: takemura@stat.t.u-tokyo.ac.jp*

We study Markov bases of decomposable graphical models consisting of primitive moves (i.e., square-free moves of degree two) by determining the structure of fibers of sample size two. We show that the number of elements of fibers of sample size two are powers of two and we characterize primitive moves in Markov bases in terms of connected components of induced subgraphs of the independence graph of a hierarchical model. This allows us to derive a complete description of minimal Markov bases and minimal invariant Markov bases for decomposable models.

## 1. Introduction

Since Sturmfels (1996) and Diaconis and Sturmfels (1998) introduced the Markov chain Monte Carlo approach based on a Markov basis for testing the goodness of fit of statistical models of multiway contingency tables, many researchers have showed the usefulness of the approach and studied Markov bases for various kinds of statistical models in computational algebraic statistics (e.g., Hoşten and Sullivant (2002); Dobra (2003); Dobra and Sullivant (2004); Geiger *et al.* (2006)). Hierarchical models are of basic importance for statistical analysis of multiway contingency tables (e.g., Lauritzen (1996); Agresti (2002)). As illustrated in Aoki and Takemura (2003), however, the structure of Markov bases for hierarchical models is very complicated in general. Decomposable models defined in terms of chordal graphs are particularly useful submodels of hierarchical models. They are known to possess Markov bases consisting of primitive moves, that is, square-free moves of degree two (Dobra (2003); Hoşten and Sullivant (2002); Geiger *et al.* (2006)). Dobra (2003) provided an algorithm to generate moves in such Markov bases based on a clique tree of the chordal graph defining the model.

The main purpose of this paper is to clarify structures of Markov bases consisting of primitive moves for decomposable models. As shown in Takemura and Aoki (2004), Markov bases for general models can be constructed by combining moves of increasing degrees. This fact indicates the importance of studying the structure of primitive moves in order to clarify the structure of Markov bases for more general hierarchical models. Some practical models such as subtable sum models (Hara *et al.* (2009)) and quasi-independence models for incomplete contingency tables that contain some structural zeros (Aoki and Takemura (2005); Rapallo (2006)) are obtained by imposing some constraints on a decomposable model. Rasch models (e.g., Chen and Small (2005)) and many-facet Rasch models (e.g., Zhu *et al.* (1998); Basturk (2008)) that are commonly used in psychometrics and behaviormetrics are considered as decomposable models restricted to contingency tables in which cell frequencies are zeros or ones. From a practical viewpoint, detailed properties of Markov bases for decomposable models may also give insights into Markov bases for such models.

The present authors have been studying Markov bases from the viewpoint of minimality (Aoki and Takemura (2003); Takemura and Aoki (2004)) and invariance (Aoki and Takemura (2008a, 2008b)) for some specific hierarchical models. The notions of minimality and invariance of Markov bases are important because they give concise expressions of the Markov basis. In this paper we extend the results to decomposable models.

The set of contingency tables sharing the same marginal frequencies corresponding to the generating set of the model is called a fiber. The structure of primitive moves is equivalent to that of fibers of sample size two. We study the structure of fibers of sample size two in detail and give a complete description of minimal Markov bases and minimal invariant Markov bases for decomposable models. We also show that construction of a minimal invariant Markov basis is directly related to a basis of a vector space over the finite field GF(2). We describe under what conditions Dobra's Markov basis is minimal or minimal invariant. We also give a necessary and sufficient condition for the uniqueness of the minimal Markov basis for decomposable models.

The organization of the paper is as follows. In Section 2 we set up notation for this paper and summarize preliminary results. In Section 3 we clarify structures of fibers of sample size two. Using this characterization, in Section 4 we give a complete description of minimal Markov bases and minimal invariant Markov bases for decomposable models. In Section 5 we briefly discuss reduced Gröbner bases for decomposable models and we end the paper with some concluding remarks in Section 6.

## 2. Preliminaries

We mostly follow the notation in Lauritzen (1996); Hoşten and Sullivant (2002); Dobra (2003) for multiway contingency tables. Let $\Delta = \{1, \ldots, m\}$ denote the set of variables of an $m$-way contingency table. Let $I_\delta$, $\delta \in \Delta$, denote the number of levels of the variable $\delta$. For convenience we take the set of levels of the variable $\delta$ as $\mathcal{I}_\delta = \{0, 1, \ldots, I_\delta - 1\}$ starting from 0 as in Hoşten and Sullivant (2002). The cells of the contingency table are indexed by

$$i = (i_1, \ldots, i_m) \in \mathcal{I} = \prod_{\delta \in \Delta} \mathcal{I}_\delta.$$

$n(i)$ denotes the frequency of the cell $i$ and $\mathbf{n} = \{n(i)\}_{i \in \mathcal{I}}$ denotes an $m$-way contingency table. The set of positive cells $\mathrm{supp}(\mathbf{n}) = \{i \in \mathcal{I} \mid n(i) > 0\}$ is the *support* of $\mathbf{n}$.

For a subset $D \subset \Delta$ of the variables, the $D$-marginal $\mathbf{n}_D$ of $\mathbf{n}$ is the contingency table with marginal cells $i_D \in \mathcal{I}_D = \prod_{\delta \in D} \mathcal{I}_\delta$ and entries given by $n_D(i_D) = \sum_{i_{D^C} \in \mathcal{I}_{D^C}} n(i_D, i_{D^C})$. Here we are denoting $i = (i_D, i_{D^C})$ by appropriately reordering indices. In this paper for notational simplicity, appropriate reordering of indices is performed as needed.

Now we consider the existence of a table $\mathbf{n}$ with the marginal tables $\mathbf{n}_{D_1}, \ldots, \mathbf{n}_{D_r}$. Dobra (2003) defined that the marginal tables $\mathbf{n}_{D_1}, \ldots, \mathbf{n}_{D_r}$ are *consistent* if, for any $r_1$, $r_2$, the $(D_{r_1} \cap D_{r_2})$-marginal of $\mathbf{n}_{D_{r_1}}$ is equal to the $(D_{r_1} \cap D_{r_2})$-marginal of $\mathbf{n}_{D_{r_2}}$. The consistency of the marginal tables is obviously a necessary condition for the existence of $\mathbf{n}$. However we note that it does not necessarily guarantee the existence of $\mathbf{n}$ in general (e.g., Irving and Jerrum (1994); Vlach (1986)).

Let $\mathcal{D} = \{D_1, \ldots, D_r\}$ be the set of facets of a simplicial complex such that $\Delta = \bigcup_{j=1}^{r} D_j$. Then $\mathcal{D}$ is called a generating class. Let $p(i)$ denote the cell probability for $i$. Then the hierarchical model for a generating class $\mathcal{D}$ is written as

$$\log p(i) = \sum_{D \in \mathcal{D}} \mu_D(i),$$

where $\mu_D$ depends only on $i_D$.

Let $\mathcal{G}^\mathcal{D}$ be a graph with the vertex set $\Delta$ and an edge between $\delta, \delta' \in \Delta$ if and only if there exists $D \in \mathcal{D}$ such that $\delta, \delta' \in D$. $\mathcal{G}^\mathcal{D}$ is called an *independence graph* of $\mathcal{D}$ (Dobra and Sullivant (2004)). A hierarchical model for $\mathcal{D}$ is called graphical if $\mathcal{D} = \{D_1, \ldots, D_r\}$ is the set of (maximal) cliques of $\mathcal{G}^\mathcal{D}$. By a clique we mean the set of vertices of a maximal complete induced subgraph. A graphical model is called *decomposable* if $\mathcal{G}^\mathcal{D}$ is chordal, that is, every cycle of $\mathcal{G}^\mathcal{D}$ with length greater than three has a chord. A *clique tree* (or a junction tree) $\mathcal{T}$ of a chordal graph $\mathcal{G}^\mathcal{D}$ is a tree, such that the vertices of $\mathcal{T}$ are cliques of $\mathcal{G}^\mathcal{D}$ and it satisfies the following property:

$$D_s \cap D_t \subset D_u \quad \text{for all } D_u \text{ on the path between } D_s \text{ and } D_t \text{ in } \mathcal{T}.$$

An intersection $S$ of neighboring cliques in a clique tree is called a minimal vertex separator. In the following $\mathcal{S}$ denotes the set of minimal vertex separators of a chordal graph. When $\mathcal{G}^\mathcal{D}$ is not connected, we regard the empty set $\varnothing$ as a minimal vertex separator of $\mathcal{G}^\mathcal{D}$.

For a clique $D \in \mathcal{D}$ of a decomposable model, let $\mathrm{Simp}(D)$ denote the set of simplicial vertices in $D$ and let $\mathrm{Sep}(D)$ denote the set of non-simplicial vertices in $D$ (Hara and Takemura (2006)). If $\text{Simp}(D) \neq \varnothing$, $D$ is called a simplicial clique. A simplicial clique $D$ is called a *boundary clique* if there exists another clique $D' \in \mathcal{D}$ such that $\text{Sep}(D) = D \cap D'$ (Shibata (1988)). Simplicial vertices in boundary cliques are called simply separated vertices (Hara and Takemura (2006)). Hara and Takemura (2006) showed that a clique $D$ is a boundary clique if and only if there exists a clique tree such that $D$ is its endpoint. Hence there exist at least two boundary cliques in any chordal graph.

Finally we summarize some relevant facts on fibers and Markov bases (Takemura and Aoki (2004, 2005)). Given the generating class $\mathcal{D} = \{D_1, \ldots, D_r\}$ of a hierarchical model, we denote the set of marginal frequencies as

$$\mathbf{b} = \{n_{D_j}(i_{D_j}), i_{D_j} \in \mathcal{I}_{D_j}, j = 1, \ldots, r\}.$$

We consider $\mathbf{b}$ as a column vector with dimension $d = \sum_{j=1}^{r} \prod_{\delta \in D_j} I_\delta$, where the elements are ordered according to an appropriate lexicographical order. We also order the elements of $\mathbf{n}$ appropriately and consider $\mathbf{n}$ as a column vector. Then the relation between the joint frequencies $\mathbf{n}$ and the marginal frequencies $\mathbf{b}$ is written simply as

$$\mathbf{b} = A\mathbf{n},$$

where $A$ is a $d \times |\mathcal{I}|$ matrix consisting of 0's and 1's. $A$ is the "incidence matrix" of cells and marginals with 1 indicating that the corresponding cell (column) is included in the corresponding marginal (row).

Given $\mathbf{b}$, the set

$$\mathcal{F}_\mathbf{b} = \{\mathbf{n} \geq 0 \mid \mathbf{b} = A\mathbf{n}\}$$

of contingency tables sharing the same marginal frequencies, $\mathbf{b}$, is called a *fiber* or $\mathbf{b}$-*fiber*, where $\mathbf{n} \geq 0$ denotes $n(i) \geq 0$ for all $i \in \mathcal{I}$. All contingency tables $\mathbf{n}$ in the same fiber $\mathcal{F}_\mathbf{b}$ have the same total frequency $n = \sum_{i \in \mathcal{I}} n(i)$. We call this common total frequency the *sample size* or the *degree* of $\mathbf{b}$ and denote it by $\deg \mathbf{b}$. We call $\mathcal{F}_\mathbf{b}$ with $\deg \mathbf{b} = 2$ a "degree two fiber" in the following.

An integer array $\mathbf{z} = \{z(i)\}_{i \in \mathcal{I}}$ of the same dimension as $\mathbf{n}$ is called a *move* if $A\mathbf{z} = \mathbf{0}$, that is, $z_D(i_D) := \sum_{i_{D^C} \in \mathcal{I}_{D^C}} z(i_D, i_{D^C}) = 0$ for all $D \in \mathcal{D}$. A move $\mathbf{z}$ is written as the difference of its positive part and negative part as $\mathbf{z} = \mathbf{z}^+ - \mathbf{z}^-$. Then $A\mathbf{z}^+ = A\mathbf{z}^-$. Therefore $\mathbf{z}^+$ and $\mathbf{z}^-$ belong to the same fiber. In this case we simply say that a move $\mathbf{z}$ belongs to the fiber $\mathcal{F}_{A\mathbf{z}^+}$. We call $\deg A\mathbf{z}^+$ the degree of a move $\mathbf{z}$. Clearly $\deg A\mathbf{z}^+ \geq 2$. Especially when $\mathbf{z}$ is a primitive move, that is, a square-free move of degree two, $\deg A\mathbf{z}^+ = 2$ and $\mathbf{z}^+$ and $\mathbf{z}^-$ belong to the same degree two fiber. Therefore the structure of primitive move is equivalent to the structure of corresponding degree two fiber. If we add a move or subtract a move $\mathbf{z}$ to $\mathbf{n} \in \mathcal{F}_\mathbf{b}$, we can move to another state $\mathbf{n} + \mathbf{z}$ (or $\mathbf{n} - \mathbf{z}$) in the same fiber $\mathcal{F}_\mathbf{b}$ as long as there is no negative element in $\mathbf{n} + \mathbf{z}$ (or $\mathbf{n} - \mathbf{z}$). A finite set $\mathcal{M}$ of moves is called a *Markov basis* if for every fiber the states become mutually accessible by the moves from $\mathcal{M}$. By using the Metropolis–Hastings procedure to control transition probabilities by moves of a Markov basis, we can construct a Markov chain on every fiber (Diaconis and Sturmfels (1998)).

A Markov basis $\mathcal{M}$ is *minimal* if every proper subset of $\mathcal{M}$ is no longer a Markov basis. Minimal Markov bases may not be unique in general. However, in view of the definition of the minimum fiber Markov basis (the set of moves that cannot be replaced by a sequence of moves of lower degree, see Section 2.2 of Takemura and Aoki (2005)), the fibers of the moves of all minimal Markov bases are common. We refer to the set of fibers common to all minimal Markov bases as the fibers of the minimum fiber Markov basis.

Suppose that a degree two fiber $\mathcal{F}_\mathbf{b}$ contains more than one element, that is, $|\mathcal{F}_\mathbf{b}| \geq 2$. Then no two elements $\mathbf{n}, \mathbf{n}'$ of the fiber share a support:

$$\deg \mathbf{b} = 2, \qquad \mathbf{n} \neq \mathbf{n}' \in \mathcal{F}_\mathbf{b} \implies \mathrm{supp}(\mathbf{n}) \cap \mathrm{supp}(\mathbf{n}') = \varnothing.$$

It follows that each element of a degree two fiber with more than one element is an indispensable monomial (Aoki *et al.* (2008)), that is, each contingency table of sample size two is isolated and has to be connected to some other table in the same fiber by a degree two move of a Markov basis. Hence each degree two fiber with more than one element has to be a fiber of the minimum fiber Markov basis. This fact holds for any hierarchical model. Note however that for some hierarchical models, such as no three-factor interaction models (Aoki and Takemura (2003)), every degree two fiber has only one element.

On the other hand, for decomposable models, Dobra (2003) has shown that there exists a Markov basis consisting of primitive moves. It implies that for decomposable models it suffices to study degree two fibers. In particular the fibers of the minimum fiber Markov bases are exactly the degree two fibers with more than one element. Furthermore, by the characterization of the uniqueness of minimal Markov bases in Takemura and Aoki (2004), it follows that the minimal Markov basis for a decomposable model is unique if and only if all degree two fibers contain at most two elements. Based on this result we will give a necessary and sufficient condition for the uniqueness of minimal Markov bases for decomposable models (Theorem 2 below) in terms of the properties of their chordal graphs.

## 3. Structure of degree two fibers

In this section we study the structure of degree two fibers to clarify the structure of primitive moves. Let $\mathcal{D} = \{D_1, \ldots, D_r\}$ be the generating class of a hierarchical model. Let $\mathbf{b}$ be a set of marginal frequencies of a contingency table with sample size two. We are interested in the structure of a degree two fiber $\mathcal{F}_\mathbf{b}$. Because the sample size is two, for each $D \in \mathcal{D}$, there exist at most two marginal cells $i_D$ with positive marginal frequency $n_D(i_D) > 0$. The same reasoning holds for each variable $\delta \in \Delta$; namely in the one-dimensional marginal table $\{n_{\{\delta\}}(i_\delta), i_\delta \in \{0, 1, \ldots, I_\delta - 1\}\}$, there exist at most two levels $i_\delta$ such that $n_{\{\delta\}}(i_\delta) > 0$. For a given $\mathbf{b}$ we say that the variable $\delta$ is *degenerate* if there exists a unique level $i_\delta$ such that $n_{\{\delta\}}(i_\delta) = 2$. Otherwise, if there exist two levels $i_\delta \neq i'_\delta$ such that $n_{\{\delta\}}(i_\delta) = n_{\{\delta\}}(i'_\delta) = 1$, then we say that the variable $\delta$ is *non-degenerate*.

If a variable $\delta$ is degenerate for a given marginal $\mathbf{b}$, then the level of the variable $\delta$ is uniquely determined from $\mathbf{b}$ and it is common for all contingency tables $\mathbf{n} \in \mathcal{F}_{\mathbf{b}}$. In particular, if all the variables $\delta \in \Delta$ are degenerate, then $\mathcal{F}_{\mathbf{b}} = \{\mathbf{n}\}$ is a one-element fiber with frequency $n(i) = 2$ at a particular cell $i$. Since this case is trivial, below we consider the case wherein at least one variable is non-degenerate. For convenience we denote

$$\mathbf{n} = (i)(j)$$

when $n(i) = n(j) = 1$, $i \neq j$. From the fact that there exist at most two levels with positive one-dimensional marginals for each variable, it follows that we only need to consider $2 \times \cdots \times 2$ tables for studying degree two fibers. Therefore, for our purposes in this section, we let $I_1 = \cdots = I_m = 2$, $\mathcal{I} = \{0,1\}^m$ without loss of generality.

For a given $\mathbf{b}$ of degree two let $\bar{\Delta}_{\mathbf{b}}$ denote the set of non-degenerate variables. As noted above we assume that $\bar{\Delta}_{\mathbf{b}} \neq \varnothing$. Each $\mathbf{n} \in \mathcal{F}_{\mathbf{b}}$ is of the form $\mathbf{n} = (i)(i') = (i_1, \ldots, i_m)(i'_1, \ldots, i'_m)$, $i \neq i'$. Furthermore, for non-degenerate $\delta \in \bar{\Delta}_{\mathbf{b}}$ the levels of the variable $\delta$ in $i$ and $i'$ are different:

$$\{i_\delta, i'_\delta\} = \{0, 1\} \qquad \forall \delta \in \bar{\Delta}_{\mathbf{b}},$$

or equivalently $i'_\delta = 1 - i_\delta$, $\forall \delta \in \bar{\Delta}_{\mathbf{b}}$. In the following we use the notation $i^*_\delta = 1 - i_\delta$. More generally for a subset $D = \{\delta_1, \ldots, \delta_k\}$ of the variables and a marginal cell $i_D = (i_{\delta_1}, \ldots, i_{\delta_k})$ we write

$$i^*_D \equiv (i^*_{\delta_1}, \ldots, i^*_{\delta_k}) = (1 - i_{\delta_1}, \ldots, 1 - i_{\delta_k}).$$

Let us identify $\mathbf{n} \in \mathcal{F}_{\mathbf{b}}$ with the set $\{i, i'\}$ of its two cells of frequency one. Then we see that the number of elements $|\mathcal{F}_{\mathbf{b}}|$ of the fiber is at most $2^{|\bar{\Delta}_{\mathbf{b}}|-1}$. However some choice of $\{i, i'\}$ with

$$i_\delta, i^*_\delta \in \{0, 1\} \qquad \forall \delta \in \bar{\Delta}_{\mathbf{b}},$$

may not be in the fiber $\mathcal{F}_{\mathbf{b}}$. This is because if $\delta$ and $\delta'$ belong to a common $D \in \mathcal{D}$, then the values of $i_\delta$ and $i_{\delta'}$ are tied together. For example, let $D = \{1, 2\} \in \mathcal{D}$ and consider the $\{1, 2\}$-marginal specified as

$$n_{\{1,2\}}(0,0) = n_{\{1,2\}}(1,1) = 1, \qquad n_{\{1,2\}}(0,1) = n_{\{1,2\}}(1,0) = 0.$$

Then if we choose $i_1 = 0$, we have to choose $i_2 = 0$. In Takemura and Hara (2007) we considered a very similar problem in the framework of swapping observations among two records in a microdata set for the purpose of statistical disclosure control. As in Takemura and Hara (2007) we make the following definition.

Let $\mathcal{G}(\bar{\Delta}_{\mathbf{b}})$ be a graph with the set of vertices $\bar{\Delta}_{\mathbf{b}}$ and an edge between $\delta \in \bar{\Delta}_{\mathbf{b}}$ and $\delta' \in \bar{\Delta}_{\mathbf{b}}$ if and only if there exists some $D \in \mathcal{D}$ such that $\delta, \delta' \in D$. Namely there exists an edge between two non-degenerate variables if and only if these two variables appear together in some marginal tables of $\mathcal{D}$. Note that $\mathcal{G}(\bar{\Delta}_{\mathbf{b}})$ is the induced subgraph of $\mathcal{G}^{\mathcal{D}}$ with the vertices restricted to $\bar{\Delta}_{\mathbf{b}}$. As discussed above in this case the values of $i_\delta$ and $i_{\delta'}$ are tied together and once the value of $i_\delta$ is chosen, for example, $i_\delta = 0$, then the value of $i_{\delta'}$ becomes fixed depending on the specifications of the marginals $n_D$.

We summarize the above argument in the following lemma.

**Lemma 1.** *Suppose that* **b** *is a set of consistent marginal frequencies of a contingency table with sample size two. Let $\Gamma$ be any subset of a connected component in $\mathcal{G}(\bar{\Delta}_\mathbf{b})$. Then the marginal table $\mathbf{n}_\Gamma = \{n_\Gamma(i_\Gamma) \,|\, i_\Gamma \in \mathcal{I}_\Gamma\}$ is uniquely determined.*

**Proof.** Let $r(\Gamma)$ be the number of generating sets $D \in \mathcal{D}$ satisfying $\Gamma \cap D \neq \varnothing$. We prove this lemma by induction on $r(\Gamma)$. When $r(\Gamma) = 1$, the lemma obviously holds. Suppose that the lemma holds for all $r(\Gamma) < r'$ and we now assume that $r(\Gamma) = r'$. Let $\Gamma_1 \subset \Gamma$ and $\Gamma_2 \subset \Gamma$ satisfy

$$\Gamma_1 \cup \Gamma_2 = \Gamma, \qquad \Gamma_1 \cap \Gamma_2 \neq \varnothing, \qquad r(\Gamma_1) < r', \qquad r(\Gamma_2) < r'.$$

Since $r(\Gamma_1) < r'$ and $r(\Gamma_2) < r'$ both $\mathbf{n}_{\Gamma_1}$ and $\mathbf{n}_{\Gamma_2}$ are uniquely determined. Suppose that

$$n_{\Gamma_1}(i_{\Gamma_1 \setminus \Gamma_2}, i_{\Gamma_1 \cap \Gamma_2}) = 1, \qquad n_{\Gamma_1}(i^*_{\Gamma_1 \setminus \Gamma_2}, i^*_{\Gamma_1 \cap \Gamma_2}) = 1. \qquad (3.1)$$

Then from the consistency of **b** there uniquely exists $i_{\Gamma_2 \setminus \Gamma_1} \in \mathcal{I}_{\Gamma_2 \setminus \Gamma_1}$ such that

$$n_{\Gamma_2}(i_{\Gamma_2 \setminus \Gamma_1}, i_{\Gamma_1 \cap \Gamma_2}) = 1, \qquad n_{\Gamma_2}(i^*_{\Gamma_2 \setminus \Gamma_1}, i^*_{\Gamma_1 \cap \Gamma_2}) = 1. \qquad (3.2)$$

Hence the table $\mathbf{n}_\Gamma = \{n(j_\Gamma) \,|\, j_\Gamma \in \mathcal{I}_\Gamma\}$ such that

$$n(j_\Gamma) = \begin{cases} 1, & \text{if } j_\Gamma = (i_{\Gamma_1 \setminus \Gamma_2}, i_{\Gamma_1 \cap \Gamma_2}, i_{\Gamma_2 \setminus \Gamma_1}) \text{ or } j_\Gamma = (i^*_{\Gamma_1 \setminus \Gamma_2}, i^*_{\Gamma_1 \cap \Gamma_2}, i^*_{\Gamma_2 \setminus \Gamma_1}), \\ 0, & \text{otherwise} \end{cases}$$

is consistent with the marginal **b**.

Suppose that there exists another marginal table $\mathbf{n}'_\Gamma$ that is consistent with **b** such that $n_\Gamma(j_\Gamma) = n_\Gamma(j^*_\Gamma) = 1$ and $j_\Gamma \neq (i_{\Gamma_1 \setminus \Gamma_2}, i_{\Gamma_1 \cap \Gamma_2}, i_{\Gamma_2 \setminus \Gamma_1})$. Then we have at least

$$n_{\Gamma_1}(i_{\Gamma_1}) = 0 \quad \text{or} \quad n_{\Gamma_2}(i_{\Gamma_2}) = 0.$$

This contradicts (3.1) and (3.2). □

By using the result of Lemma 1, we obtain the following theorem on the number of elements in degree two fibers.

**Theorem 1.** *Let $\mathcal{F}_\mathbf{b}$ be a degree two fiber such that $\bar{\Delta}_\mathbf{b} \neq \varnothing$ and let $c(\mathbf{b})$ be the number of connected components of $\mathcal{G}(\bar{\Delta}_\mathbf{b})$. Then*

$$|\mathcal{F}_\mathbf{b}| = 2^{c(\mathbf{b})-1}.$$

**Proof.** Denote by $\Gamma_1, \ldots, \Gamma_c$, $c = c(\mathbf{b})$, the connected components of $\mathcal{G}(\bar{\Delta}_\mathbf{b})$. Define $\Gamma_{c+1}$ by $\Gamma_{c+1} = \Delta \setminus \bar{\Delta}_\mathbf{b}$. Then there exists $i_{\Gamma_{c+1}}$ such that

$$i_{\Gamma_{c+1}} = \{i_\delta \mid \delta \in \Gamma_{c+1}, n_{\{\delta\}}(i_\delta) = 2\}.$$

From Lemma 1 the marginal cells $i_{\Gamma_k}$ such that $n_{\Gamma_k}(i_{\Gamma_k}) = n_{\Gamma_k}(i^*_{\Gamma_k}) = 1$ uniquely exist for $k = 1, \ldots, c$. Now define $\mathcal{I}_\mathbf{b}$ by

$$\mathcal{I}_\mathbf{b} = \{i_{\Gamma_1}, i^*_{\Gamma_1}\} \times \{i_{\Gamma_2}, i^*_{\Gamma_2}\} \times \cdots \times \{i_{\Gamma_c}, i^*_{\Gamma_c}\} \times \{i_{\Gamma_{c+1}}\},$$

where $\times$ denotes the direct product of sets. Suppose that $j \in \mathcal{I}_\mathbf{b}$. Define $\mathbf{n}_j = \{n_j(i) \mid i \in \mathcal{I}\}$ by

$$n_j(i) = \begin{cases} 1, & \text{if } i = j \text{ or } i = j^*, \\ 0, & \text{otherwise.} \end{cases}$$

Then we have $\mathcal{F}(\mathcal{I}_\mathbf{b}) = \{\mathbf{n}_j \mid j \in \mathcal{I}_\mathbf{b}\} \subset \mathcal{F}_\mathbf{b}$ and $|\mathcal{F}(\mathcal{I}_\mathbf{b})| = 2^{c-1}$.

If there exists $\mathbf{n}' = \{n'(i) \mid i \in \mathcal{I}\}$ such that $\mathbf{n}' \in \mathcal{F}_\mathbf{b}$ and $\mathbf{n}' \notin \mathcal{F}(\mathcal{I}_\mathbf{b})$, then there exists a cell $j \in \mathcal{I}$ and $1 \leq k \leq c+1$ such that $n(j) = 1$ and $j_{\Gamma_k} \neq i_{\Gamma_k}$. This implies that there exists $D_l \in \mathcal{D}$ such that $n'(i_{D_l}) \neq n(i_{D_l})$. Hence we have $|\mathcal{F}_\mathbf{b}| = 2^{c(\mathbf{b})-1}$. □

As mentioned in Section 2, for a consistent $\mathbf{b}$ such that $\deg \mathbf{b} > 2$, it is known that $\mathcal{F}_\mathbf{b}$ may be empty (e.g., Irving and Jerrum (1994); Vlach (1986)) in general. However Theorem 1 shows that, in the case $\deg \mathbf{b} = 2$, if a consistent $\mathbf{b}$ such that $\bar{\Delta}_\mathbf{b} \neq \varnothing$ is given, then $\mathcal{F}_\mathbf{b} \neq \varnothing$ for any hierarchical model.

It is helpful to consider permuting the levels $0 \leftrightarrow 1$ for each variable and understand Theorem 1 in a canonical form. This amounts to considering invariance of hierarchical models with respect to permutation of levels of each variable as studied in Aoki and Takemura (2008a). Although we have reduced our consideration to $2^m$ tables in treating degree two fibers, we are really considering general hierarchical models of $I_1 \times \cdots \times I_m$ tables. Note that hierarchical models possess the symmetry with respect to relabeling the levels of each variable, that is, it is invariant under the action of the direct product of symmetric groups $S_{I_1} \times \cdots \times S_{I_m}$ acting on the set of cells. If we again restrict our consideration to degree two fibers, we only need to consider the action of $S_2^m = S_2 \times \cdots \times S_2$. It is clear that structures of degree two fibers are invariant under the action of $S_2^m$.

In particular as a "representative fiber", we can consider $\mathbf{b}$ such that the levels of all degenerate variables are determined as 0. Also for such a $\mathbf{b}$, let $\Gamma \subset \bar{\Delta}_\mathbf{b}$ be the set of vertices of a connected component of $\mathcal{G}(\bar{\Delta}_\mathbf{b})$. Then we can without loss of generality assume that two $\Gamma$-marginal cells of frequency 1 are specified as

$$1 = n_\Gamma(0, 0, \ldots, 0) = n_\Gamma(1, 1, \ldots, 1). \tag{3.3}$$

This can be achieved by interchanging the levels of each variable in $\bar{\Delta}_\mathbf{b}$. Under this standardization the proof of Theorem 1 is easier to understand, because for each connected component of $\mathcal{G}(\bar{\Delta}_\mathbf{b})$ we either choose all 0's or all 1's for the component.

This standardization is also useful in determining the setwise stabilizer of $\mathcal{F}_\mathbf{b}$ in $S_2^m$ (Section 3.1 of Aoki and Takemura (2008b)). If we standardize the levels as (3.3), then the setwise stabilizer of $\mathcal{F}_\mathbf{b}$ is isomorphic to $c(\mathbf{b})$-fold direct product of $S_2$'s:

$$S_2^{c(\mathbf{b})} = S_2 \times \cdots \times S_2.$$

In other words the structure of $\mathcal{F}_\mathbf{b}$ is equivalent to the structure of the fiber $\mathcal{F}_{\mathbf{b}'}$ with $\Delta = \bar\Delta_{\mathbf{b}'} = \{1, \ldots, c(\mathbf{b}')\}$ and totally disconnected $\mathcal{G}(\bar\Delta_{\mathbf{b}'})$. In the next section we use this fact in determining the minimal invariant Markov bases for decomposable models.

Finally we prove the following theorem on a sufficient condition for non-uniqueness of minimal Markov bases.

**Theorem 2.** *Let $\mathcal{D} = \{D_1, \ldots, D_r\}$ be the generating class of a hierarchical model. Suppose that $m \geq 3$ and there exist three variables $\delta_1, \delta_2, \delta_3$ that are not connected to each other in $\mathcal{G}^\mathcal{D}$. Then minimal Markov bases for the hierarchical model with the generating class $\mathcal{D}$ are not unique.*

**Proof.** It suffices to find a degree two fiber with more than two elements. Consider $\mathbf{b}$ such that $\bar\Delta_\mathbf{b} = \{\delta_1, \delta_2, \delta_3\}$. From the condition of the theorem $\mathcal{G}(\bar\Delta_\mathbf{b})$ has an induced subgraph with three connected components. Therefore $|\mathcal{F}_\mathbf{b}| = 4$. This completes the proof. □

## 4. Markov bases for decomposable models

### 4.1. Minimal and unique minimal Markov bases

In this section we discuss Markov bases of decomposable models in detail from the viewpoint of minimality based on the results obtained in the previous section. Since there exists a Markov basis consisting of primitive moves for decomposable models, the set of fibers of the minimum fiber Markov basis coincides with the set of degree two fibers with more than one element. Theorem 1 of the previous section enables a complete description of minimal Markov bases of decomposable models.

Let $\deg \mathbf{b} = 2$. Let $\mathcal{T}_\mathbf{b}$ be any tree whose nodes are elements of $\mathcal{F}_\mathbf{b}$. Denote the set of edges in $\mathcal{T}_\mathbf{b}$ by $\mathcal{M}_{\mathcal{T}_\mathbf{b}}$. We note that we can identify each edge $(\mathbf{n}, \mathbf{n}') \in \mathcal{M}_{\mathcal{T}_\mathbf{b}}$ with a move $\mathbf{z} = \mathbf{n} - \mathbf{n}'$. So we identity $\mathcal{M}_{\mathcal{T}_\mathbf{b}}$ with a set of moves for $\mathcal{F}_\mathbf{b}$. In considering Markov bases, we ignore the sign of $\mathbf{z}$ and identify $\mathbf{z} = \mathbf{n} - \mathbf{n}'$ with $-\mathbf{z} = \mathbf{n}' - \mathbf{n}$ and consider the edges in $\mathcal{T}_\mathbf{b}$ as undirected. In contrast when we consider Gröbner bases, we distinguish $\mathbf{z}$ from $-\mathbf{z}$ and correspondingly consider directed edges.

Let $\mathcal{B}_{\mathrm{nd}}$ be

$$\mathcal{B}_{\mathrm{nd}} = \{\mathbf{b} \mid \deg \mathbf{b} = 2, |\mathcal{F}_\mathbf{b}| \geq 2\}. \tag{4.1}$$

Then we define $\mathcal{M}^0$ as follows,

$$\mathcal{M}^0 = \bigcup_{\mathbf{b} \in \mathcal{B}_{\mathrm{nd}}} \mathcal{M}_{\mathcal{T}_\mathbf{b}}. \tag{4.2}$$

By following Dobra (2003) and Takemura and Aoki (2004), we easily obtain the following theorem.

**Theorem 3.** *$\mathcal{M}^0$ is a minimal Markov basis and (4.2) is a disjoint union. Conversely every minimal Markov basis can be written as in (4.2).*

***Example 1 (The complete independence model of three-way contingency tables).*** Consider the model $\mathcal{D} = \{\{1\}, \{2\}, \{3\}\}$ for the $2 \times 2 \times 2$ contingency tables. $\mathcal{B}_{\mathrm{nd}}$ for the model has seven elements. Denote them by $\mathbf{b}_1, \ldots, \mathbf{b}_7$. Figure 1 shows an example of $\mathcal{M}_{\mathcal{T}_{\mathbf{b}_t}}$ for $t = 1, \ldots, 7$. $\mathbf{b}_1, \ldots, \mathbf{b}_7$ satisfy

$$\bar{\Delta}_{\mathbf{b}_1} = \{1, 2, 3\}, \qquad \bar{\Delta}_{\mathbf{b}_2} = \bar{\Delta}_{\mathbf{b}_3} = \{1, 2\},$$
$$\bar{\Delta}_{\mathbf{b}_4} = \bar{\Delta}_{\mathbf{b}_5} = \{2, 3\}, \qquad \bar{\Delta}_{\mathbf{b}_6} = \bar{\Delta}_{\mathbf{b}_7} = \{1, 3\}. \tag{4.3}$$

The union of all these moves is a minimal Markov basis for the model. Since $\mathcal{F}_{\mathbf{b}_1}$ is a four elements fiber, $\mathcal{T}_{\mathbf{b}_1}$ is not uniquely determined. Hence minimal Markov bases are not unique for this model.

As seen from this example, minimal Markov bases are not necessarily uniquely determined. Based on Theorems 1 and 3, we can derive a necessary and sufficient condition on decomposable models to have the unique minimal Markov basis.

**Corollary 1.** *There exists the unique minimal Markov basis for a decomposable model if and only if the number of connected components in any induced subgraphs of $\mathcal{G}^{\mathcal{D}}$ is less than three.*

**Proof.** Suppose that $\mathcal{G}(\bar{\Delta}_{\mathbf{b}})$ has more than two connected components. Then since $|\mathcal{F}_{\mathbf{b}}| \geq 4$ from Theorem 1, $\mathcal{T}_{\mathbf{b}}$ is not uniquely determined. For a different tree $\mathcal{T}'_{\mathbf{b}}$, $\mathcal{M}_{\mathcal{T}_{\mathbf{b}}} \neq \mathcal{M}_{\mathcal{T}'_{\mathbf{b}}}$. Hence minimal Markov bases are not unique either.

Conversely, assume that the number of connected components of $\mathcal{G}(\bar{\Delta}_{\mathbf{b}})$ for all $\mathbf{b} \in \mathcal{B}_{\mathrm{nd}}$ is two. Then $\mathcal{T}_{\mathbf{b}}$ for all $\mathbf{b} \in \mathcal{B}_{\mathrm{nd}}$ is uniquely determined. Hence the minimal Markov basis is unique. □

For decomposable models $\mathcal{G}^{\mathcal{D}}$ is chordal. From the graph theoretical viewpoint the above corollary can be rewritten as follows.

**Corollary 2.** *For a decomposable model, there exists the unique minimal Markov basis if and only if $\mathcal{G}^{\mathcal{D}}$ has only two boundary cliques $D$ and $D'$ such that $D'' \subset D \cup D'$ for all $D'' \in \mathcal{D}$.*

**Proof.** Suppose that $\mathcal{G}^{\mathcal{D}}$ has two boundary cliques $D$ and $D'$ such that $D'' \subset D \cup D'$ for all $D'' \in \mathcal{D}$. Then any vertex in $D''$ is adjacent to $D$ or $D'$. Hence the number of connected components for any induced subgraph of $\mathcal{G}^{\mathcal{D}}$ is at most two.

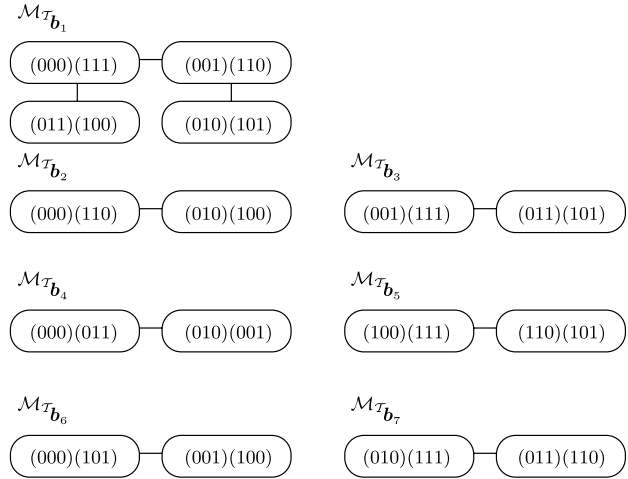

(The triplets in brackets refer to cells in the contingency table.)

**Figure 1.** $\mathcal{M}_{\mathcal{T}_{\mathbf{b}_t}}$ in the complete independence model of three-way contingency tables.

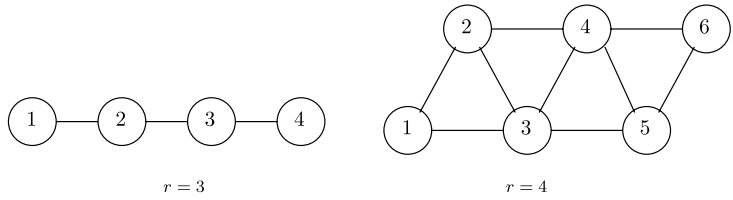

**Figure 2.** Examples of the graphs satisfying the condition of Theorem 2.

Conversely suppose that there exists $D'' \in \mathcal{D}$ such that $D'' \not\subset D \cup D'$. Then the subgraph induced by the union of $D'' \setminus (D \cup D')$, $\mathrm{Simp}(D)$ and $\mathrm{Simp}(D')$ has three connected components. □

The graphs with $r = 2$ always satisfy the conditions of the theorem. For $r \geq 3$ the graph with

$$\mathcal{D} = \{\{1,\ldots,r-1\},\{2,\ldots,r\},\ldots,\{r,\ldots,2r-2\}\} \qquad (4.4)$$

satisfies the conditions of the theorem. Figure 2 represents the graphs satisfying (4.4) for $r = 3, 4$. We can easily see that any induced subgraph of the graphs in the figure has at most two connected components.

Let $\mathcal{T} = (\mathcal{D}, \mathcal{E})$ be a clique tree for $\mathcal{G}^{\mathcal{D}}$. Denote by $\mathcal{T}_e = (\mathcal{D}_e, \mathcal{E}_e)$ and $\mathcal{T}'_e = (\mathcal{D}'_e, \mathcal{E}'_e)$ the two induced subtrees of $\mathcal{T}$ obtained by removing an edge $e \in \mathcal{E}$. Define $V_e$ and $V'_e$ by

$$V_e = \bigcup_{D \in \mathcal{D}_e} D, \qquad V'_e = \bigcup_{D \in \mathcal{D}'_e} D.$$

Let $\mathcal{M}^{\mathcal{T}}(V_e, V'_e)$ be the set of all primitive moves for the decomposable model determined by the chordal graph whose set of cliques is $\{V_e, V'_e\}$. Dobra (2003) showed that

$$\mathcal{M}^{\mathcal{T}} = \bigcup_{e \in \mathcal{E}} \mathcal{M}^{\mathcal{T}}(V_e, V'_e) \tag{4.5}$$

is a Markov basis. We call $\mathcal{M}^{\mathcal{T}}$ a Dobra's Markov basis. From the viewpoint of minimality of Markov bases, we have the following theorem.

**Theorem 4.** *A decomposable model has a clique tree $\mathcal{T}$ such that $\mathcal{M}^{\mathcal{T}}$ is a minimal Markov basis if and only if the model has the unique minimal Markov basis.*

**Proof.** When a decomposable model has unique minimal Markov basis, $\mathcal{M}^{\mathcal{T}}$ coincides with the minimal Markov basis.

Suppose that there exist three vertices in $\mathcal{G}^{\mathcal{D}}$ that are not adjacent to each other. Let $1, 2$ and $3$ be such three vertices and assume that $l \in D_l$, $D_l \in \mathcal{D}$, for $l = 1, 2, 3$. Define $\{1,2,3\}^c = \Delta \setminus \{1,2,3\}$. Consider a degree two fiber $\mathcal{F}_{\mathbf{b}}$ such that $\bar{\Delta}_{\mathbf{b}} = \{1,2,3\}$ and $n_{\{1,2,3\}^c}(i_{\{1,2,3\}^c}) = 2$ for some $i_{\{1,2,3\}^c}$. Then $|\mathcal{F}_{\mathbf{b}}| = 4$ from Theorem 1 and we can denote these four elements by

$$\begin{aligned}
\mathbf{n}_1 &= (000 i_{\{1,2,3\}^c})(111 i_{\{1,2,3\}^c}), \\
\mathbf{n}_2 &= (001 i_{\{1,2,3\}^c})(110 i_{\{1,2,3\}^c}), \\
\mathbf{n}_3 &= (010 i_{\{1,2,3\}^c})(101 i_{\{1,2,3\}^c}), \\
\mathbf{n}_4 &= (011 i_{\{1,2,3\}^c})(100 i_{\{1,2,3\}^c}).
\end{aligned} \tag{4.6}$$

A minimal Markov basis connects these four elements by three moves. Let $\mathcal{T} = (\mathcal{D}, \mathcal{E})$ be any clique tree for $\mathcal{G}^{\mathcal{D}}$ and $\mathcal{T}' = (\mathcal{D}', \mathcal{E}')$ be the smallest subtree of $\mathcal{T}$ satisfying $D_l \in \mathcal{D}'$ for $l = 1, 2$ and $3$. Then we can assume that $\mathcal{T}'$ satisfies either of the following two conditions,

(i) $D_2$ is an interior point and $D_1$ and $D_3$ are endpoints on the path;
(ii) all of $D_1, D_2$ and $D_3$ are endpoints of $\mathcal{T}'$.

In both cases there exists $e \in \mathcal{E}$ such that $D_1, D_2 \subset V_e$ and $D_3 \subset V'_e$. Then $\mathcal{M}^{\mathcal{T}}(V_e, V'_e)$ includes the following two moves,

$$\mathbf{z}_1 = \mathbf{n}_1 - \mathbf{n}_2, \qquad \mathbf{z}_2 = \mathbf{n}_3 - \mathbf{n}_4.$$

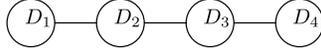

**Figure 3.** $\mathcal{T}$ in Example 2.

On the other hand, there also exists $e' \in \mathcal{E}$ such that $D_1 \subset V_{e'}$ and $D_2, D_3 \subset V_{e'}'$. In this case $\mathcal{M}^{\mathcal{T}}(V_{e'}, V_{e'}')$ includes the following two moves,

$$\mathbf{z}_3 = \mathbf{n}_1 - \mathbf{n}_4, \qquad \mathbf{z}_4 = \mathbf{n}_2 - \mathbf{n}_3.$$

Thus $\mathcal{M}^{\mathcal{T}}$ includes at least four moves for the fiber $\mathcal{F}_\mathbf{b}$, which implies that $\mathcal{M}^{\mathcal{T}}$ is not minimal for the model that does not have the unique minimal Markov basis. □

***Example 2 (The complete independence model of four-way contingency tables).*** Consider the $2 \times 2 \times 2 \times 2$ complete independence model $\mathcal{D} = \{\{1\}, \{2\}, \{3\}, \{4\}\}$. Let $\mathcal{F}_\mathbf{b}$ be the fiber with $\bar{\Delta}_\mathbf{b} = \{1, 2, 3, 4\}$, that is, $c(\mathbf{b}) = 4$ and $|\mathcal{F}_\mathbf{b}| = 8$. Consider $\mathcal{M}^{\mathcal{T}}$ for $\mathcal{T}$ in Figure 3. Denote the set of moves for $\mathcal{F}_\mathbf{b}$ belonging to $\mathcal{M}^{\mathcal{T}}$ by $\mathcal{M}^{\mathcal{T}}_\mathbf{b}$. Figure 4 shows $\mathcal{M}^{\mathcal{T}}_\mathbf{b}$. As seen from Figure 4, $\mathcal{M}^{\mathcal{T}}_\mathbf{b}$ includes 12 moves. Since $|\mathcal{F}_\mathbf{b}| = 8$, 7 moves are sufficient to connect $\mathcal{F}_\mathbf{b}$.

### 4.2. Minimal invariant Markov bases

In this section we consider Markov bases from the viewpoint of invariance under the action of the product of symmetric groups $G = G_{I_1,\ldots,I_m} = S_{I_1} \times \cdots \times S_{I_m}$ on the levels of the variables. The organization of this section is as follows. We first express a minimal invariant Markov basis as a union of orbits of $G_{I_1,\ldots,I_m}$, which minimally connects representative fibers (see (4.7) below). Then we show that the minimal set of orbits connecting a non-degenerate fiber is in one-to-one correspondence to a basis of a vector space over the finite field GF(2) (Lemma 2 and Theorem 5 below). Then the structure of minimal invariant Markov bases is given in Theorem 6.

According to Aoki and Takemura (2008a), a set of moves $\mathcal{M}$ is called *G-invariant* if

$$g \in G, \quad \mathbf{z} \in \mathcal{M} \quad \Rightarrow \quad g(\mathbf{z}) \in \mathcal{M} \quad \text{or} \quad -g(\mathbf{z}) \in \mathcal{M}.$$

$\mathcal{M}$ is called a *G-invariant Markov basis* for $\mathcal{D}$ if it is a Markov basis and also *G*-invariant. A *G*-invariant Markov basis $\mathcal{M}$ is *minimal* if no proper *G*-invariant subset of $\mathcal{M}$ is a Markov basis.

As discussed at the end of Section 3, by appropriate reordering of the indices we can consider a representative fiber

$$\mathcal{F}^0_\mathbf{b} \ni \mathbf{n}^\mathbf{b}_0 \equiv (0\cdots 0)(1\cdots 1 \overbrace{0\cdots 0}^{|\Delta \setminus \bar{\Delta}_\mathbf{b}|}).$$

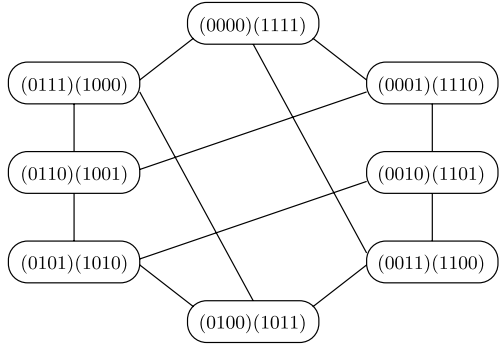

**Figure 4.** $\mathcal{M}_{\mathbf{b}}^{\mathcal{T}}$ for $\mathbf{b}$ such that $\bar{\Delta}_{\mathbf{b}} = \{1,2,3,4\}$.

Then any $\mathbf{n} \in \mathcal{F}_{\mathbf{b}}^0$ is expressed as follows,

$$\mathbf{n} = (\overbrace{0\cdots 0}^{|\Gamma_1|} i_{\Gamma_2} \cdots i_{\Gamma_{c(\mathbf{b})}} \overbrace{0\cdots 0}^{|\Delta\setminus\bar{\Delta}_{\mathbf{b}}|})(\overbrace{1\cdots 1}^{|\Gamma_1|} i^*_{\Gamma_2} \cdots i^*_{\Gamma_{c(\mathbf{b})}} \overbrace{0\cdots 0}^{|\Delta\setminus\bar{\Delta}_{\mathbf{b}}|}),$$

$$i_{\Gamma_l} = \overbrace{0\cdots 0}^{|\Gamma_l|} \quad \text{or} \quad i_{\Gamma_l} = \overbrace{1\cdots 1}^{|\Gamma_l|}, \qquad l = 2,\ldots,c(\mathbf{b}),$$

where $\Gamma_l$ are the connected components of $\mathcal{G}(\bar{\Delta}_{\mathbf{b}})$. Let $G^{\Gamma_l}$, $l=2,\ldots,c(\mathbf{b})$, be the diagonal subgroup of $S_2^{|\Gamma_l|}$ defined by

$$G^{\Gamma_l} = \{\bar{g} = (g,\ldots,g) \mid g \in S_2\} \subset S_2^{|\Gamma_l|}.$$

Define

$$G_{\mathbf{b}} = G^{\Gamma_2} \times \cdots \times G^{\Gamma_{c(\mathbf{b})}}$$

and let $g \in G_{\mathbf{b}}$ act on $\mathbf{n} \in \mathcal{F}_{\mathbf{b}}^0$ by

$$g(\mathbf{n}) = (\overbrace{0\cdots 0}^{|\Gamma_1|} \bar{g}_2(i_{\Gamma_2}) \cdots \bar{g}_{c(\mathbf{b})}(i_{\Gamma_{c(\mathbf{b})}}) \overbrace{0\cdots 0}^{|\Delta\setminus\bar{\Delta}_{\mathbf{b}}|})(\overbrace{1\cdots 1}^{|\Gamma_1|} \bar{g}_2(i^*_{\Gamma_2}) \cdots \bar{g}_{c(\mathbf{b})}(i^*_{\Gamma_{c(\mathbf{b})}}) \overbrace{0\cdots 0}^{|\Delta\setminus\bar{\Delta}_{\mathbf{b}}|}).$$

Clearly $g(\mathbf{n}) \in \mathcal{F}_{\mathbf{b}}^0$ for $\mathbf{n} \in \mathcal{F}_{\mathbf{b}}^0$ and furthermore for any $\mathbf{n} \in \mathcal{F}_{\mathbf{b}}^0$ there exists $g \in G_{\mathbf{b}}$ such that $\mathbf{n} = g(\mathbf{n}_0^{\mathbf{b}})$. This shows that $G_{\mathbf{b}} \subset G_{I_1,\ldots,I_m}$ is the setwise stabilizer of $\mathcal{F}_{\mathbf{b}}^0$ acting transitively on $\mathcal{F}_{\mathbf{b}}^0$.

Let $\mathcal{M}_{G_{\mathbf{b}}}$ be a minimal $G_{\mathbf{b}}$-invariant set of moves that connects $\mathcal{F}_{\mathbf{b}}^0$. Let $\kappa(\mathbf{b})$ be the number of $G_{\mathbf{b}}$-orbits included in $\mathcal{M}_{G_{\mathbf{b}}}$. As representative moves of $G_{\mathbf{b}}$-orbits in $\mathcal{M}_{G_{\mathbf{b}}}$ we can consider

$$\mathbf{z}_k^{\mathbf{b}} = \mathbf{n}_0^{\mathbf{b}} - \mathbf{n}_k^{\mathbf{b}} \in \mathcal{M}_{\mathbf{b}}, \qquad \mathbf{n}_k^{\mathbf{b}} \in \mathcal{F}_{\mathbf{b}}^0, \qquad k = 1,\ldots,\kappa(\mathbf{b}).$$

This is because we can always send $\mathbf{n}$ in $\mathbf{z} = \mathbf{n} - \mathbf{n}'$ to $\mathbf{n}_0^{\mathbf{b}}$ by the transitivity of $G_{\mathbf{b}}$. Denote $\mathcal{M}_{G_{\mathbf{b}}}^0 = \{\mathbf{z}_1^{\mathbf{b}}, \ldots, \mathbf{z}_{\kappa(\mathbf{b})}^{\mathbf{b}}\}$. Define the set of representative fibers by

$$\mathcal{B}_{\mathrm{nd}}^0 = \{\mathbf{b} \mid \mathbf{n}_0^{\mathbf{b}} \in \mathcal{F}_{\mathbf{b}}^0\} \subset \mathcal{B}_{\mathrm{nd}}.$$

From Aoki and Takemura (2008a) a minimal $G_{I_1,\ldots,I_m}$-invariant Markov basis can be expressed by

$$\mathcal{M}_G = \bigcup_{\mathbf{b} \in \mathcal{B}_{\mathrm{nd}}^0} \bigcup_{k=1}^{\kappa(\mathbf{b})} G_{I_1,\ldots,I_m}(\mathbf{z}_k^{\mathbf{b}}), \tag{4.7}$$

where $G_{I_1,\ldots,I_m}(\mathbf{z}_k^{\mathbf{b}})$ denotes the $G_{I_1,\ldots,I_m}$-orbit through $\mathbf{z}_k^{\mathbf{b}}$. Hence in order to clarify the structure of $\mathcal{M}_G$, it suffices to consider $2 \times \cdots \times 2$ tables and investigate $\kappa(\mathbf{b})$ and $\mathcal{M}_{G_{\mathbf{b}}}^0$ for each $\mathcal{F}_{\mathbf{b}}^0$.

As mentioned in Section 3, the structure of $\mathcal{F}_{\mathbf{b}}^0$ is equivalent to the structure of the fiber with $\bar{\Delta}_{\mathbf{b}} = \Delta = \{1, \ldots, c(\mathbf{b})\}$ and $\mathcal{G}(\bar{\Delta}_{\mathbf{b}})$ is totally disconnected. We first consider the structure of such a fiber. $\mathcal{F}_{\mathbf{b}}^0$ satisfies

$$\mathcal{F}_{\mathbf{b}}^0 = \{(0i_2 \cdots i_{c(\mathbf{b})})(1i_2^* \cdots i_{c(\mathbf{b})}^*) \mid (i_2 \cdots i_{c(\mathbf{b})}) = i_{\Delta \setminus \{1\}} \in \mathcal{I}_{\Delta \setminus \{1\}}\} \tag{4.8}$$

and $(0 \cdots 0)(1 \cdots 1) \in \mathcal{F}_{\mathbf{b}}^0$. Then we note that we can identify $G_{\mathbf{b}}$ with $S_2^{c(\mathbf{b})-1}$. For $g \in S_2^{c(\mathbf{b})-1}$, we write $g = (g_2, \ldots, g_{c(\mathbf{b})})$, where $g_l \in S_2$ for $l = 2, \ldots, c(\mathbf{b})$. A representative move of $S_2^{c(\mathbf{b})-1}$-orbit is expressed by

$$\mathbf{z}^{\mathbf{b}} = (0 \cdots 0)(1 \cdots 1) - (0i_{\Delta \setminus \{1\}})(1i_{\Delta \setminus \{1\}}^*)$$

for some $i_{\Delta \setminus \{1\}} \in \mathcal{I}_{\Delta \setminus \{1\}}$. We first derive $\kappa(\mathbf{b})$ and $\mathcal{M}_{G_{\mathbf{b}}}$. Let $\mathcal{V}^{c(\mathbf{b})-1} = \{0, 1\}^{c(\mathbf{b})-1}$ denote the $(c(\mathbf{b}) - 1)$-dimensional vector space over the finite field GF(2), where the addition of two vectors is defined as the "exclusive or" (XOR) of the elements. Let $\oplus$ denote the XOR operation. Let $\circ$ denote the group operation of $S_2^{c(\mathbf{b})-1}$. Then we obtain the following lemma.

**Lemma 2.** $S_2^{c(\mathbf{b})-1}$ *is isomorphic to* $\mathcal{V}^{c(\mathbf{b})-1}$.

**Proof.** Consider the map $\phi: S_2^{c(\mathbf{b})-1} \to \mathcal{V}^{c(\mathbf{b})-1}$ such that $\phi(g) = \mathbf{v} = (v_2, \ldots, v_{c(\mathbf{b})}) \in \mathcal{V}^{c(\mathbf{b})-1}$, where

$$v_l = \begin{cases} 0, & \text{if } g_l(i_l) = i_l, \\ 1, & \text{if } g_l(i_l) = i_l^* \end{cases}$$

for $l = 2, \ldots, c(\mathbf{b})$ and $\{i_l, i_l^*\} = \{0, 1\}$. For $g' = (g_2', \ldots, g_{c(\mathbf{b})}') \in S_2^{c(\mathbf{b})-1}$, $g_l' \in S_2$, and $\mathbf{v}' \in \mathcal{V}^{c(\mathbf{b})-1}$, define $\phi(g') = \mathbf{v}' = (v_2', \ldots, v_{c(\mathbf{b})}')$. Then we have $\phi(g \circ g') = \tilde{\mathbf{v}} = (\tilde{v}_2, \ldots, \tilde{v}_{c(\mathbf{b})})$,

$\tilde{\mathbf{v}} \in \mathcal{V}^{c(\mathbf{b})-1}$, where

$$\tilde{v}_l = \begin{cases} 0, & \text{if } g_l \circ g_l'(i_l) = i_l, \\ 1, & \text{if } g_l \circ g_l'(i_l) = i_l^* \end{cases}$$

for $l = 2, \ldots, c(\mathbf{b})$. Hence we have

$$\tilde{v}_l = v_l \oplus v_l', \qquad l = 2, \ldots, c(\mathbf{b}).$$

Therefore $\phi$ is a homomorphism. It is obvious that $\phi$ is a bijection. Therefore $S_2^{c(\mathbf{b})-1}$ is isomorphic to $\mathcal{V}^{c(\mathbf{b})-1}$. $\square$

Based on this lemma, we can show the equivalence between $S_2^{c(\mathbf{b})-1}$-orbits in a minimal $S_2^{c(\mathbf{b})-1}$-invariant set of moves that connects $\mathcal{F}_\mathbf{b}^0$ and a (vector space) basis of $\mathcal{V}^{c(\mathbf{b})-1}$.

**Theorem 5.** *Let* $\mathcal{V}^0 = \{\mathbf{v}_k = (v_{k2}, \ldots, v_{kc(\mathbf{b})}), k = 2, \ldots, c(\mathbf{b})\}$ *be any basis of* $\mathcal{V}^{c(\mathbf{b})-1}$. *Define* $\mathbf{n}_0^\mathbf{b}, \mathbf{n}_{\mathbf{v}_k}^\mathbf{b} \in \mathcal{F}_\mathbf{b}^0$ *by*

$$\mathbf{n}_0^\mathbf{b} = (00 \cdots 0)(11 \cdots 1), \qquad \mathbf{n}_{\mathbf{v}_k}^\mathbf{b} = (0 v_{k2} \cdots v_{kc(\mathbf{b})})(1 v_{k2}^* \cdots v_{kc(\mathbf{b})}^*),$$

*where* $v_{kl}^* = 1 \oplus v_{kl}$. *Let* $\mathcal{M}_{G_\mathbf{b}}$ *be an* $S_2^{c(\mathbf{b})-1}$-*invariant set of moves in* $\mathcal{F}_\mathbf{b}^0$. *Then* $\mathcal{M}_{G_\mathbf{b}}$ *is a minimal* $S_2^{c(\mathbf{b})-1}$-*invariant set of moves that connects* $\mathcal{F}_\mathbf{b}^0$ *if and only if the representative moves of the* $S_2^{c(\mathbf{b})-1}$-*orbits in* $\mathcal{M}_{G_\mathbf{b}}$ *are expressed by* $\mathbf{z}_{\mathbf{v}_k}^\mathbf{b} = \mathbf{n}_0^\mathbf{b} - \mathbf{n}_{\mathbf{v}_k}^\mathbf{b}$, $k = 2, \ldots, c(\mathbf{b})$. *Hence* $\kappa(\mathbf{b}) = c(\mathbf{b}) - 1$.

**Proof.** Suppose that $\mathcal{M}_{G_\mathbf{b}}$ is a minimal $S_2^{c(\mathbf{b})-1}$-invariant set of moves that connects $\mathcal{F}_\mathbf{b}$ and that $\mathcal{M}_{G_\mathbf{b}}$ includes $\kappa(\mathbf{b})$ orbits as $S_2^{c(\mathbf{b})-1}(\mathbf{z}_1^\mathbf{b}), \ldots, S_2^{c(\mathbf{b})-1}(\mathbf{z}_{\kappa(\mathbf{b})}^\mathbf{b})$, where

$$\mathbf{z}_k^\mathbf{b} = \mathbf{n}_0^\mathbf{b} - \mathbf{n}_k^\mathbf{b}, \qquad \mathbf{n}_k^\mathbf{b} = (0 i_{k2} \cdots i_{kc(\mathbf{b})})(1 i_{k2}^* \cdots i_{kc(\mathbf{b})}^*)$$

for $i_{kl} \in \mathcal{I}_l$, $k = 1, \ldots, \kappa(\mathbf{b})$, $l = 2, \ldots, c(\mathbf{b})$. Let $g^k \in S_2^{c(\mathbf{b})-1}$ satisfy $g^k(\mathbf{n}_0^\mathbf{b}) = \mathbf{n}_k^\mathbf{b}$ for $k = 1, \ldots, \kappa(\mathbf{b})$. We write $g^k = (g_{k2}, \ldots, g_{kc(\mathbf{b})})$, $g_{kl} \in S_2$ for $l = 2, \ldots, c(\mathbf{b})$. Let $H_\mathbf{b} = \{g^1, \ldots, g^{\kappa(\mathbf{b})}\} \subset S_2^{c(\mathbf{b})-1}$ be a subset of $S_2^{c(\mathbf{b})-1}$. As mentioned above, $\mathcal{F}_\mathbf{b}^0$ can be expressed as in (4.8). Hence for any $\mathbf{n} \in \mathcal{F}_\mathbf{b}^0$ there exists $g \in S_2^{c(\mathbf{b})-1}$ satisfying $\mathbf{n} = g(\mathbf{n}_0^\mathbf{b})$. $\mathcal{M}_{G_\mathbf{b}}$ connects $\mathcal{F}_\mathbf{b}^0$ if and only if there exists $p \leq \kappa(\mathbf{b})$ such that

$$\mathbf{n} = \mathbf{n}_0^\mathbf{b} - \mathbf{z}_{k_1}^\mathbf{b} - g^{k_1}(\mathbf{z}_{k_2}^\mathbf{b}) - \cdots - g^{k_{p-1}} \circ \cdots \circ g^{k_1}(\mathbf{z}_{k_p}^\mathbf{b})$$

and $g = g^{k_p} \circ \cdots \circ g^{k_1}$. Hence $\mathcal{M}_{G_\mathbf{b}}$ is a minimal $S_2^{c(\mathbf{b})-1}$-invariant set of moves that connects $\mathcal{F}_\mathbf{b}$ if and only if $H_\mathbf{b}$ satisfies

$$\forall g \in S_2^{c(\mathbf{b})-1}, \exists p \leq \kappa(\mathbf{b}), \exists g^{k_1} \in H_\mathbf{b}, \ldots, \exists g^{k_p} \in H_\mathbf{b} \qquad \text{s.t. } g = g^{k_p} \circ \cdots \circ g^{k_1} \qquad (4.9)$$

and no proper subset of $H_\mathbf{b}$ satisfies (4.9).

Denote $\mathcal{V}^0 = \phi(H_\mathbf{b}) \subset \mathcal{V}^{c(\mathbf{b})-1}$. Then from Lemma 2, (4.9) is equivalent to

$$\forall \mathbf{v} \in \mathcal{V}, \exists \mathbf{v}_2 \in \mathcal{V}^0, \ldots, \exists \mathbf{v}_{p+1} \in \mathcal{V}^0 \qquad \text{s.t. } \mathbf{v} = \mathbf{v}_2 \oplus \cdots \oplus \mathbf{v}_{p+1}. \tag{4.10}$$

From the minimality of $\mathcal{M}_{G_\mathbf{b}}$ no proper subset of $\mathcal{V}^0$ satisfies (4.10). This implies that $\mathcal{V}^0$ is a basis of $\mathcal{V}^{c(\mathbf{b})-1}$ and hence $\kappa(\mathbf{b}) = c(\mathbf{b}) - 1$. If we define $g^k = \phi^{-1}(\mathbf{v}_{k+1})$ for $k = 1, \ldots, c(\mathbf{b}) - 1$, we have $g_{kl}(0) = v_{k+1,l}$ and hence $g^k(\mathbf{n}_0^\mathbf{b}) = \mathbf{n}_k^\mathbf{b} = \mathbf{n}_{\mathbf{v}_{k+1}}^\mathbf{b}$. Therefore $\mathbf{z}_{\mathbf{v}_k}^\mathbf{b}$, $k = 2, \ldots, c(\mathbf{b})$, are the representative moves of the $S_2^{c(\mathbf{b})-1}$-orbits in $\mathcal{M}_{G_\mathbf{b}}$.

Conversely suppose that the representative moves of $\mathcal{M}_{G_\mathbf{b}}$ are $\mathbf{z}_{\mathbf{v}_k}^\mathbf{b}$, $k = 2, \ldots, c(\mathbf{b})$. $\mathcal{V}^0$ satisfies (4.10) and no proper subset of $\mathcal{V}^0$ satisfies (4.10). Hence if we define $g^k = \phi^{-1}(\mathbf{v}_{k+1})$ and $H_\mathbf{b} = \{g^1, \ldots, g^{c(\mathbf{b})-1}\}$, $H_\mathbf{b}$ satisfies (4.9) and no proper subset of $H_\mathbf{b}$ satisfies (4.9). Hence $\mathcal{M}_{G_\mathbf{b}}$ is a minimal $S_2^{c(\mathbf{b})-1}$-invariant set of moves that connects $\mathcal{F}_\mathbf{b}$. □

For example, we can set $\mathcal{V}^0 = \{\mathbf{v}_2, \ldots, \mathbf{v}_{c(\mathbf{b})}\}$ as

$$\mathbf{v}_2 = (11\cdots11), \qquad \mathbf{v}_3 = (01\cdots11), \qquad \ldots,$$
$$\mathbf{v}_{c(\mathbf{b})-1} = (00\cdots011), \qquad \mathbf{v}_{c(\mathbf{b})} = (00\cdots01),$$

and then the representative moves in a minimal $G$-invariant Markov basis are

$$\begin{aligned}
\mathbf{z}_2^0 &= (00\cdots0)(11\cdots1) - (011\cdots11)(100\cdots00), \\
\mathbf{z}_3^0 &= (00\cdots0)(11\cdots1) - (001\cdots11)(110\cdots00), \\
&\vdots \qquad \vdots \qquad \vdots \\
\mathbf{z}_{c(\mathbf{b})}^0 &= (00\cdots0)(11\cdots1) - (000\cdots01)(111\cdots10).
\end{aligned} \tag{4.11}$$

So far we have focused on $\mathcal{F}_\mathbf{b}$ such that $\bar{\Delta}_\mathbf{b} = \Delta = \{1, \ldots, c(\mathbf{b})\}$ and $\mathcal{G}(\bar{\Delta}_\mathbf{b})$ is totally disconnected. Now we consider a fiber for a general $\mathbf{b}$ of a general decomposable model. Define $\bar{g}_{kl} \in G^{\Gamma_l}$ by

$$\bar{g}_{kl}(\overbrace{0\cdots0}^{|\Gamma_l|}) = \begin{cases} 0\cdots0, & \text{if } v_{k+1} = 0, \\ 1\cdots1, & \text{if } v_{k+1} = 1 \end{cases} \tag{4.12}$$

for $k = 1, \ldots, c(\mathbf{b}) - 1$ and $l = 2, \ldots, c(\mathbf{b})$ and define $g^k \in G_\mathbf{b}$ by

$$\begin{aligned}
g^k(\mathbf{n}) = (\overbrace{0\cdots0}^{|\Gamma_1|} \bar{g}_{k2}(i_{\Gamma_2}) \cdots \bar{g}_{kc(\mathbf{b})}(i_{\Gamma_{c(\mathbf{b})}}) \overbrace{0\cdots0}^{|\Delta\setminus\bar{\Delta}_\mathbf{b}|}) \\
\times (\overbrace{1\cdots1}^{|\Gamma_1|} \bar{g}_{k2}(i^*_{\Gamma_2}) \cdots \bar{g}_{kc(\mathbf{b})}(i^*_{\Gamma_{c(\mathbf{b})}}) \overbrace{0\cdots0}^{|\Delta\setminus\bar{\Delta}_\mathbf{b}|}).
\end{aligned} \tag{4.13}$$

Denote $\mathbf{n}_{\mathbf{v}_{k+1}}^\mathbf{b} = g^k(\mathbf{n}_0^\mathbf{b})$ and $\mathbf{z}_{\mathbf{v}_{k+1}}^\mathbf{b} = \mathbf{n}_0^\mathbf{b} - \mathbf{n}_{\mathbf{v}_{k+1}}^\mathbf{b}$. Based on (4.7) and Theorem 5, we can easily obtain the following result.

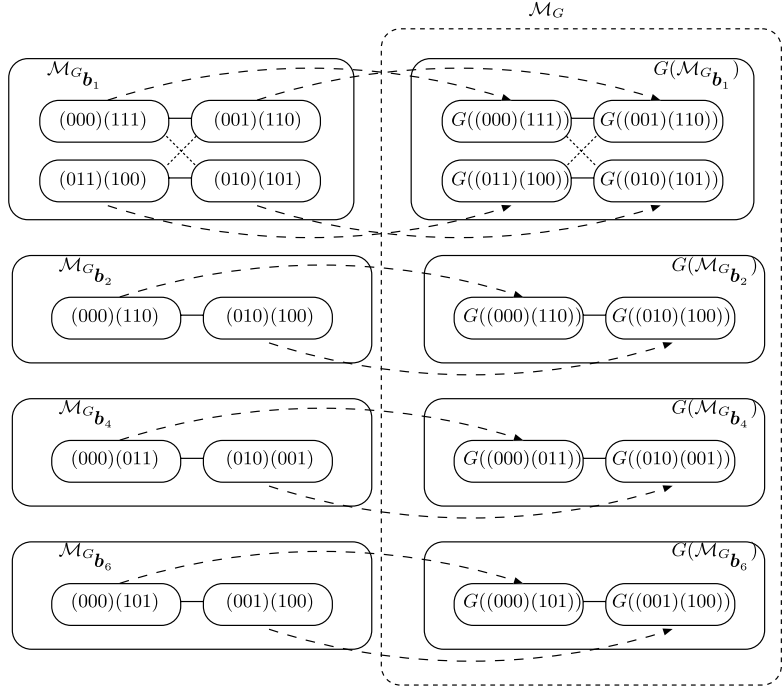

**Figure 5.** The structure of minimal $G_{2,2,2}$-invariant Markov bases for the complete independence model of three-way contingency tables.

**Theorem 6.** $\mathcal{M}_{G_\mathbf{b}}$ *is a minimal $S_2^{c(\mathbf{b})-1}$-invariant set of moves that connects $\mathcal{F}_\mathbf{b}^0$ if and only if the representative moves of the $S_2^{c(\mathbf{b})-1}$-orbits in $\mathcal{M}_{G_\mathbf{b}}$ are expressed as $\mathbf{z}_{\mathbf{v}_k}^\mathbf{b}$, $k = 2, \ldots, c(\mathbf{b})$. Hence $\kappa(\mathbf{b}) = c(\mathbf{b}) - 1$. Then*

$$\mathcal{M}_G = \bigcup_{\mathbf{b} \in \mathcal{B}_{\mathrm{nd}}^0} \bigcup_{k=2}^{c(\mathbf{b})} G_{I_1,\ldots,I_m}(\mathbf{z}_k^\mathbf{b})$$

*is a minimal $G_{I_1,\ldots,I_m}$-invariant Markov basis. Conversely, every minimal $G_{I_1,\ldots,I_m}$-invariant Markov basis can be written in this form.*

**Example 3 (The complete independence model of three-way contingency tables).** Define $\mathbf{b}_t$ as in Figure 1 of Example 1. Then $\mathcal{B}_{\mathrm{nd}}^0 = \{\mathbf{b}_1, \mathbf{b}_2, \mathbf{b}_4, \mathbf{b}_6\}$. Figure 5 shows a structure of $\mathcal{M}_G$ for the complete independence model of $2 \times 2 \times 2$ contingency tables. The left half of the figure shows the structure of $\mathcal{M}_{G_{\mathbf{b}_t}}$, $\mathbf{b}_t \in \mathcal{B}_{\mathrm{nd}}^0$.

$c(\mathbf{b}_1) = 3$ and hence $\kappa(\mathbf{b}_1) = 2$. If we set $\mathbf{v}_1^\mathbf{b} = (10)$ and $\mathbf{v}_2^\mathbf{b} = (01)$, we have

$$\mathbf{z}_1^{\mathbf{b}_1} = (000)(111) - (010)(101), \qquad \mathbf{z}_2^{\mathbf{b}_1} = (000)(111) - (001)(110).$$

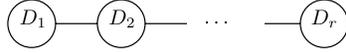

**Figure 6.** The clique tree with two endpoints.

The orbits $S_2^2(\mathbf{z}_1^{\mathbf{b}_1})$ and $S_2^2(\mathbf{z}_2^{\mathbf{b}_1})$ are expressed in dotted lines and solid lines, respectively, in the figure.

$c(\mathbf{b}_t) = 2$ and $\kappa(\mathbf{b}_t) = 1$ for $t = 2, 4, 6$. There exists one orbit in $\mathcal{M}_{G_{\mathbf{b}_t}}$ for $t = 2, 4, 6$. Then from Theorem 6 a minimal $G_{2,2,2}$-invariant Markov basis is expressed by

$$\mathcal{M}_G = G(\mathbf{z}_1^{\mathbf{b}_1}) \cup G(\mathbf{z}_2^{\mathbf{b}_1}) \cup G(\mathbf{z}_1^{\mathbf{b}_2}) \cup G(\mathbf{z}_1^{\mathbf{b}_4}) \cup G(\mathbf{z}_1^{\mathbf{b}_6}).$$

Next we consider a Dobra's Markov basis $\mathcal{M}^\mathcal{T}$ from the viewpoint of invariance. Since $\mathcal{M}^\mathcal{T}$ does not depend on the levels of the variables, $\mathcal{M}^\mathcal{T}$ is $G_{I_1,\ldots,I_m}$-invariant. Based on the result of Theorem 5, we can show that $\mathcal{M}^\mathcal{T}$ is not always a minimal invariant Markov basis.

**Theorem 7.** $\mathcal{M}^\mathcal{T}$ *is minimal invariant if and only if* $\mathcal{T}$ *has only two endpoints.*

**Proof.** It suffices to show that the theorem holds for $2 \times \cdots \times 2$ tables. Suppose that $\mathcal{T} = (\mathcal{D}, \mathcal{E})$ has more than two endpoints. Let $D_1$, $D_2$ and $D_3$ be three of them. Then they are boundary cliques. Suppose $1, 2, 3 \in \Delta$ are simply separated vertices in $D_1$, $D_2$ and $D_3$, respectively. In the same way as the argument in the proof of Theorem 4, there exist $e, e', e'' \in \mathcal{E}$ such that

$$D_1, D_2 \in V_e, \qquad D_3 \in V_e',$$
$$D_2, D_3 \in V_{e'}, \qquad D_1 \in V_{e'}',$$
$$D_3, D_1 \in V_{e''}, \qquad D_2 \in V_{e''}'.$$

Consider the moves for the fiber $\mathcal{F}_\mathbf{b}^0$ for $\mathbf{b}$ such that $\bar{\Delta}_\mathbf{b} = \{1, 2, 3\}$. Define $\mathbf{z}_5$ and $\mathbf{z}_6$ by

$$\mathbf{z}_5 = \mathbf{n}_1 - \mathbf{n}_3, \qquad \mathbf{z}_6 = \mathbf{n}_2 - \mathbf{n}_4,$$

where $\mathbf{n}_1, \ldots, \mathbf{n}_4$ are defined in (4.6). Then we have

$$\mathbf{z}_1, \mathbf{z}_2 \in \mathcal{M}^\mathcal{T}(V_e, V_e'), \qquad \mathbf{z}_3, \mathbf{z}_4 \in \mathcal{M}^\mathcal{T}(V_{e'}, V_{e'}'), \qquad \mathbf{z}_5, \mathbf{z}_6 \in \mathcal{M}^\mathcal{T}(V_{e''}, V_{e''}').$$

We note that $\{\mathbf{z}_1, \mathbf{z}_2\}$, $\{\mathbf{z}_3, \mathbf{z}_4\}$ and $\{\mathbf{z}_5, \mathbf{z}_6\}$ are $S_2^2$-orbits in $\mathcal{M}_\mathbf{b}^\mathcal{T}$. Since $\kappa(\mathbf{b}) = 2$, $\mathcal{M}^\mathcal{T}$ is not minimal invariant.

Suppose that $\mathcal{T}$ has only two endpoints. Then $\mathcal{T}$ is expressed as in Figure 6. Let $\Gamma_1(\mathbf{b}), \ldots, \Gamma_{c(\mathbf{b})}(\mathbf{b})$ be the $c(\mathbf{b})$ connected components of $\mathcal{G}(\bar{\Delta}_\mathbf{b})$. Suppose $\delta_l \in \Gamma_l(\mathbf{b})$. The structure of $\mathcal{F}_\mathbf{b}^0$ is equivalent to the structure of $\mathcal{F}_{\mathbf{b}'}^0$ such that $\bar{\Delta}_{\mathbf{b}'} = \{\delta_1, \ldots, \delta_{c(\mathbf{b})-1}\}$ and $\mathcal{G}(\bar{\Delta}_{\mathbf{b}'})$ is totally disconnected. So we restrict our consideration to such a fiber. Denote by $\mathcal{F}^0_{\mathbf{b}'}$ the representative fiber for $\mathbf{b}'$. Let

$$\mathcal{M}_{\mathbf{b}'} = \{\mathbf{n} - \mathbf{n}' \mid \mathbf{n}, \mathbf{n}' \in \mathcal{F}^0_{\mathbf{b}'}, \mathbf{n} \neq \mathbf{n}'\}$$

denote the set of all moves in $\mathcal{F}^0_{\mathbf{b}'}$. Without loss of generality we can assume that $\delta_l \in D_{\pi(l)}$, where $\pi(1) < \cdots < \pi(c(\mathbf{b}'))$. Define $e_l = (D_{l-1}, D_l) \in \mathcal{E}$, $S_l = D_{l-1} \cap D_l$, $V_l = V_{e_l} \setminus S_l$ and $V'_l = V'_{e_l} \setminus S_l$ for $l = 2, \ldots, c(\mathbf{b}')$. Then the moves in $\mathcal{M}^{\mathcal{T}}(V_l, V'_l)$ are expressed by

$$\mathbf{z} = (i_{V_l}, i_{V'_l}, i_{S_l})(j_{V_l}, j_{V'_l}, i_{S_l}) - (i_{V_l}, j_{V'_l}, i_{S_l})(j_{V_l}, i_{V'_l}, i_{S_l}),$$
$$i_{V_l}, j_{V_l} \in \mathcal{I}_{V_l}, i_{V'_l}, j_{V'_l} \in \mathcal{I}_{V'_l}, i_{S_l} \in \mathcal{I}_{S_l}.$$

If $V_{e_l} \cap \bar{\Delta}_{\mathbf{b}'} = \varnothing$ or $V'_{e_l} \cap \bar{\Delta}_{\mathbf{b}'} = \varnothing$, then we have $\mathcal{M}^{\mathcal{T}}(V_{e_l}, V'_{e_l}) \cap \mathcal{M}_{\mathbf{b}'} = \varnothing$. If $V_{e_l} \cap \bar{\Delta}_{\mathbf{b}'} \neq \varnothing$ and $V'_{e_l} \cap \bar{\Delta}_{\mathbf{b}'} \neq \varnothing$, then there exists $2 \leq k(e_l) \leq c(\mathbf{b}')$ satisfying $\delta_k \in V_l$ for all $k < k(e_l)$ and $\delta_k \in V'_l$ for all $k \geq k(e_l)$. Then

$$\mathcal{M}^{\mathcal{T}}(V_{e_l}, V'_{e_l}) \cap \mathcal{M}_{\mathbf{b}'} = S_2^{c(\mathbf{b})-1}(\mathbf{z}^0_{k(e_l)}),$$

where $\mathbf{z}^0_{k(e_l)}$ is defined as in (4.11). Hence we have

$$\mathcal{M}^{\mathcal{T}}_{\mathbf{b}'} = \bigcup_{e_l \in \mathcal{E}} \mathcal{M}^{\mathcal{T}}(V_{e_l}, V'_{e_l}) \cap \mathcal{M}_{\mathbf{b}'} = \bigcup_{k=2}^{c(\mathbf{b}')} S_2^{c(\mathbf{b})-1}(\mathbf{z}^0_k),$$

which includes $c(\mathbf{b}') - 1$ orbits for all $\mathbf{b}' \in \mathcal{B}^0_{\mathrm{nd}}$. Hence $\mathcal{M}^{\mathcal{T}}$ is minimal $G_{I_1, \ldots, I_m}$-invariant. □

*Example 4 (The complete independence model of four-way contingency tables).* As an example we consider the $2 \times 2 \times 2 \times 2$ complete independence model $\mathcal{D} = \{D_l = \{l\}, l = 1, \ldots, 4\}$. Both $\mathcal{T}^1$ and $\mathcal{T}^2$ in Figure 7 are clique trees for $\mathcal{D}$. From Theorem 7, $\mathcal{M}^{\mathcal{T}^1}$ is a minimal $S_2^3$-invariant Markov basis. Consider the representative fiber $\mathcal{F}^0_{\mathbf{b}}$ such that $\bar{\Delta}_{\mathbf{b}} = \{1, 2, 3\}$. For $j = 1, 2$, denote the two induced subtrees of $\mathcal{T}^j$ obtained by removing the edge $e_l$ by $\mathcal{T}^j_{e_l}$ and $\mathcal{T}^{j'}_{e_l}$. Figure 8 shows $\mathcal{T}^1_{e_l}$, $\mathcal{T}^{1'}_{e_l}$ and $\mathcal{M}^{\mathcal{T}^1}(V_{e_l}, V'_{e_l}) \cap \mathcal{M}_{\mathbf{b}}$. If we remove $e_3$ from $\mathcal{T}^1$, 1, 2 and 3 are still connected and hence $\mathcal{M}^{\mathcal{T}^1}(V_{e_3}, V'_{e_3}) \cap \mathcal{M}_{\mathbf{b}} = \varnothing$. Therefore $\mathcal{M}^{\mathcal{T}^1}_{\mathbf{b}}$ includes $\kappa(\mathbf{b}) = 2$ orbits.

On the other hand, since $\mathcal{T}^2$ has three endpoints, $\mathcal{M}^{\mathcal{T}^2}$ is not a minimal $S_2^3$-invariant Markov basis. Figure 9 shows $\mathcal{T}^2_{e_l}$, $\mathcal{T}^{2'}_{e_l}$ and $\mathcal{M}^{\mathcal{T}^2}(V_{e_l}, V'_{e_l}) \cap \mathcal{M}_{\mathbf{b}}$. We can see that $\mathcal{M}^{\mathcal{T}^2}_{\mathbf{b}}$ includes three orbits. As seen from this example, in general the minimality of $\mathcal{M}^{\mathcal{T}}$ depends on clique trees $\mathcal{T}$.

*Example 5.* We consider the model defined by the chordal graph in Figure 10. The clique tree of this graph is uniquely determined by $\mathcal{T}^2$ in Figure 7. As seen from this example, there exist decomposable models such that $\mathcal{M}^{\mathcal{T}}$ for every clique tree $\mathcal{T}$ is not minimal $G_{I_1, \ldots, I_m}$-invariant.

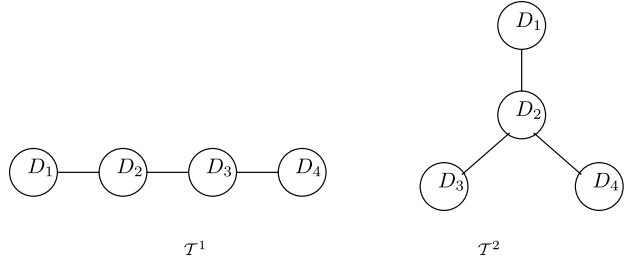

**Figure 7.** Clique trees for the 4-way complete independence model.

### 4.3. The relation between minimal and minimal invariant Markov bases

From a practical point of view a $G_{I_1,\ldots,I_m}$-invariant Markov basis is useful because its representative moves give the most concise expression of a Markov basis. On the other hand, a minimal Markov basis is also important because the number of moves contained in it is minimum among Markov bases. Here we consider the relation between a minimal and a minimal $G_{I_1,\ldots,I_m}$-invariant Markov basis and give an algorithm to obtain a minimal Markov basis from representative moves of a minimal $G_{I_1,\ldots,I_m}$-invariant Markov basis.

As mentioned in the previous section, the set of $G_{\mathbf{b}}$-orbits in a minimal $G_{\mathbf{b}}$-invariant set, $\mathcal{M}_{G_{\mathbf{b}}}$, of moves that connects $\mathcal{F}_{\mathbf{b}}^0$ has one-to-one correspondence to a basis $\mathcal{V}^0$ of $\mathcal{V}^{c(\mathbf{b})-1}$. Define $\bar{g}_{kl} \in G_{\Gamma_l}$ and $g^k \in G_\mathbf{b}$ as in (4.12) and (4.13). Let $H_\mathbf{b} = \{g^1, \ldots, g^{c(\mathbf{b})-1}\} \subset G_\mathbf{b}$. Now we generate a set of moves $\mathcal{M}_\mathbf{b}^*$ in $\mathcal{F}_\mathbf{b}$ by the following algorithm.

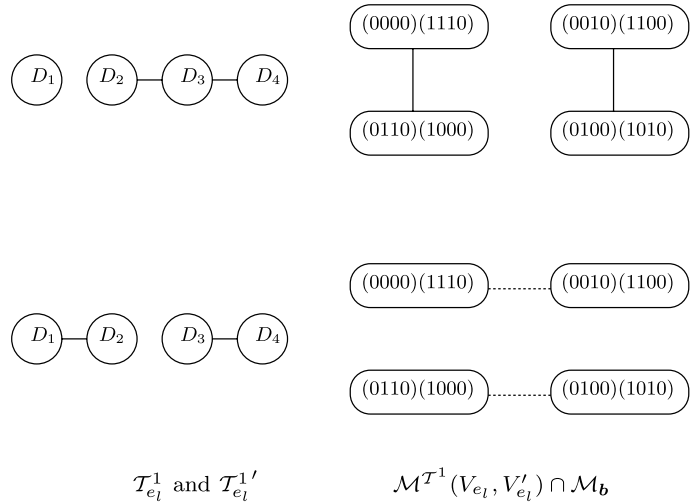

**Figure 8.** The structure of $\mathcal{M}_{\mathbf{b}}^{\mathcal{T}^1}$.

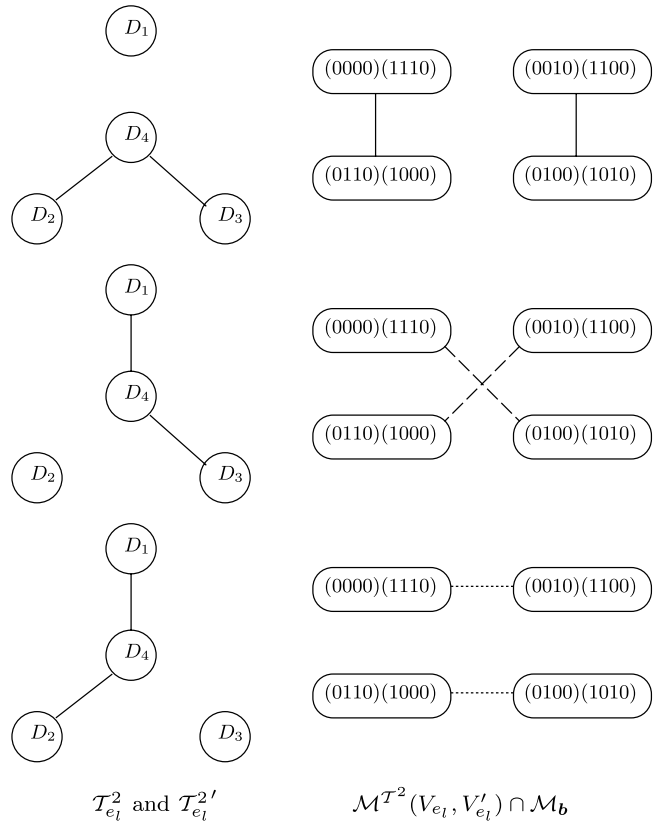

Figure 9. The structure of $\mathcal{M}_{\mathbf{b}}^{\mathcal{T}^2}$.

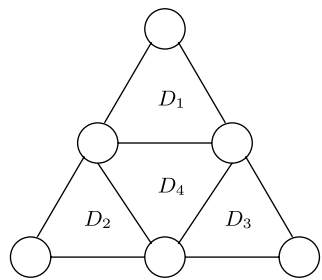

Figure 10. A chordal graph whose clique tree is uniquely determined.

*Algorithm 1.*
Input: $\mathcal{F}_\mathbf{b}$, $H_\mathbf{b} = \{g^1, \ldots, g^{c(\mathbf{b})-1}\}$
Output: $\mathcal{M}_\mathbf{b}^*$

begin
    $\mathcal{M}_\mathbf{b}^* \leftarrow \varnothing$;
    Choose any element $\mathbf{n}_1$ in $\mathcal{F}_\mathbf{b}$;
    for $k = 2$ to $c(\mathbf{b})$ do
    begin
        for $l = 1$ to $2^{k-2}$ do
        begin
            $\mathbf{n}_{l+2^{k-2}} := g^{k-1}(\mathbf{n}_l)$;
            $\mathbf{z}_{l+2^{k-2}} := \mathbf{n}_l - \mathbf{n}_{l+2^{k-2}}$;
            $\mathcal{M}_\mathbf{b}^* \leftarrow \mathcal{M}_\mathbf{b}^* \cup \{\mathbf{z}_{l+2^{k-2}}\}$;
        end
    end
    return $\mathcal{M}_\mathbf{b}^*$;
end.

**Theorem 8.** $\mathcal{M}_\mathbf{b}^*$ *generated by Algorithm 1 is a minimal set of moves that connects* $\mathcal{F}_\mathbf{b}$.

**Proof.** Since $|\mathcal{M}_\mathbf{b}^*| = 2^0 + 2^1 + \cdots + 2^{c(\mathbf{b})-1} = 2^{c(\mathbf{b})-1} - 1$, it suffices to show that $\mathbf{n}_l \neq \mathbf{n}_{l'}$ for $l \neq l'$. Suppose that there exist $l$ and $l'$, $l \neq l'$, such that $\mathbf{n}_l = \mathbf{n}_{l'}$ and $\mathbf{n}_l$, $\mathbf{n}_{l'}$ are expressed by

$$\mathbf{n}_l = g^{k_p} \circ \cdots \circ g^{k_1}(\mathbf{n}_1), \qquad \mathbf{n}_{l'} = g^{k'_{p'}} \circ \cdots \circ g^{k'_1}(\mathbf{n}_1),$$

where $k_1 < k_2 < \cdots < k_p \leq c(\mathbf{b}) - 1$ and $k'_1 < k'_2 < \cdots < k'_{p'} \leq c(\mathbf{b}) - 1$. Without loss of generality we can assume $p < p'$. Then we have

$$g^{k_p} \circ \cdots \circ g^{k_1} = g^{k'_{p'}} \circ \cdots \circ g^{k'_1} \qquad (4.14)$$

and there exists $l \leq p$ such that $k_l \neq k'_l$. From Lemma 2 (4.14) is equivalent to

$$\mathbf{v}_{k_1} \oplus \cdots \oplus \mathbf{v}_{k_p} = \mathbf{v}_{k'_1} \oplus \cdots \oplus \mathbf{v}_{k'_{p'}},$$

which contradicts that $\mathcal{V}^0$ is a basis of $\mathcal{V}^{c(\mathbf{b})-1}$. Hence we have $\mathbf{n}_l \neq \mathbf{n}_{l'}$ for $l \neq l'$. $\square$

From (4.2) we obtain the following result.

**Corollary 3.** $\mathcal{M}^* = \bigcup_{\mathbf{b} \in \mathcal{B}_{\mathrm{nd}}} \mathcal{M}_\mathbf{b}^*$ *is a minimal Markov basis.*

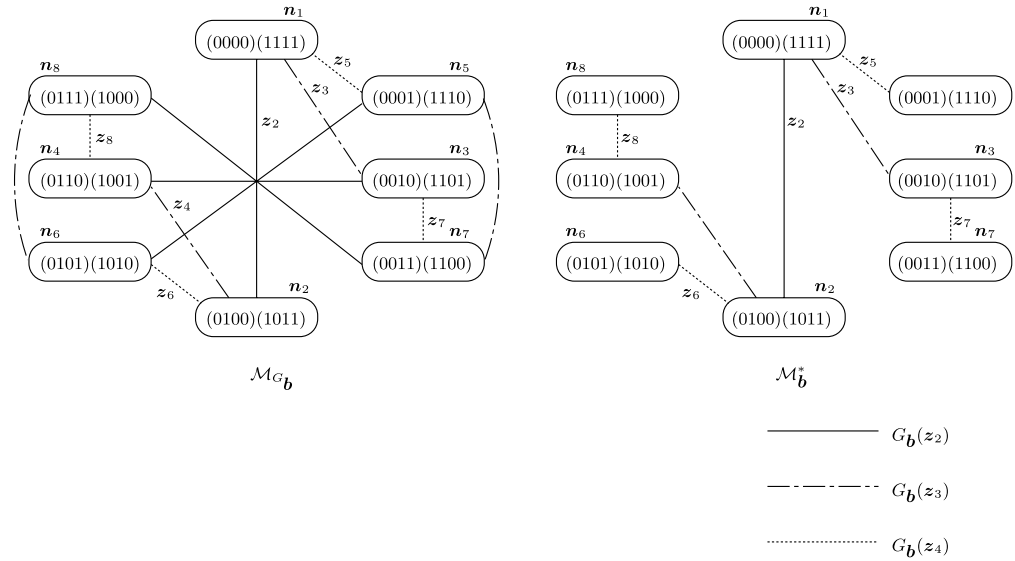

**Figure 11.** $\mathcal{M}_{G_{\mathbf{b}}}$ and $\mathcal{M}_{\mathbf{b}}^*$ generated by Algorithm 1.

*Example 6 (The complete independence model of a four-way contingency table).* We consider the same fiber as in Example 2. Define $\mathcal{V}^0 = \{\mathbf{v}_2, \mathbf{v}_3, \mathbf{v}_4\}$ by $\mathbf{v}_2 = (100)$, $\mathbf{v}_3 = (010)$ and $\mathbf{v}_4 = (001)$. Figure 11 shows $\mathcal{M}_{G_{\mathbf{b}}}$ and $\mathcal{M}_{\mathbf{b}}^*$ generated by Algorithm 1 with $\mathbf{n}_1 = (0000)(1111)$.

## 5. Gröbner bases for decomposable models

So far we have been discussing Markov bases. In this section we briefly discuss Gröbner bases. For decomposable models, Theorem 4.17 of Hoşten and Sullivant (2002) gives a recursive method for determining the term order and the corresponding Gröbner basis consisting of primitive moves only. It gives a Gröbner basis version of Dobra's Markov basis in (4.5). In Theorem 4 we saw that Dobra's construction gives a minimal Markov basis only in a special case. The same phenomenon can be observed with respect to the *reducedness* of the Gröbner basis if we simply apply Theorem 4.17 of Hoşten and Sullivant recursively, that is, the operation of Theorem 4.17 of Hoşten and Sullivant does not preserve reducedness in general. Here we are interested in an explicit description of appropriate term order and the reduced Gröbner basis for decomposable models. We prove that for decomposable models, there exists a term order such that the reduced Gröbner basis is explicitly described and, furthermore, it is minimal as a Markov basis.

In obtaining a nice Gröbner basis, the term order has to be carefully chosen. For example, consider the simple case of $3 \times 3$ two-way contingency tables with fixed row sums and column sums. Proposition 5.4 of Sturmfels (1996) shows that the set of nine primitive moves of the form

$$\pm \begin{array}{|cc|} \hline +1 & -1 \\ -1 & +1 \\ \hline \end{array}$$

constitute a reduced Gröbner basis when the cells are lexicographically ordered and the term order is chosen to be the reverse lexicographic term order. However, if we order the nine cells as

| 1 | 8 | 6 |
|---|---|---|
| 4 | 2 | 9 |
| 7 | 5 | 3 |

and use the lexicographic order, then the reduced Gröbner contains the following degree three move

$$\begin{array}{|ccc|} \hline 0 & -1 & +1 \\ +1 & 0 & -1 \\ -1 & +1 & 0 \\ \hline \end{array}$$

in addition to the nine primitive moves. This example shows that the existence of a reduced Gröbner basis consisting of primitive moves depends on the choice of a term order.

We need several steps in constructing a nice term order for a decomposable model of an $m$-way contingency table. First, we order $m$ variables. Choose a boundary clique of the chordal graph corresponding to the decomposable model and order the variables in the boundary cliques as the lowest variables. Then remove the boundary clique from the chordal graph, choose a boundary clique from the smaller graph and order the variables from the boundary clique as the next lowest variables. By recursively removing boundary cliques we obtain an ordering of variables. The resulting order is a perfect elimination scheme but has a stronger property. Second, given the order of the variables, we order the cells of an $m$-way contingency table lexicographically. Finally, as the term order $\succ$ we use the reverse lexicographic term order.

Let $\mathcal{B}_{\mathrm{nd}}$ as in (4.1). In each fiber $\mathcal{F}_{\mathbf{b}}$, $\mathbf{b} \in \mathcal{B}_{\mathrm{nd}}$, there exists the lowest element $\mathbf{n}_{\mathbf{b}}^*$ with respect to the above term order $\succ$. Define

$$\mathcal{M}^{\mathrm{GB}} = \bigcup_{\mathbf{b} \in \mathcal{B}_{\mathrm{nd}}} \bigcup_{\substack{\mathbf{n} \in \mathcal{F}_{\mathbf{b}} \\ \mathbf{n} \neq \mathbf{n}_{\mathbf{b}}^*}} \{\mathbf{n} - \mathbf{n}_{\mathbf{b}}^*\}.$$

Then we have the following theorem.

**Theorem 9.** *$\mathcal{M}^{\mathrm{GB}}$ is the reduced Gröbner basis and it is minimal as a Markov basis.*

We omit the details of the proof. By generalizing the proof of Proposition 5.4 of Sturmfels (1996) we can show that $\mathcal{M}^{\mathrm{GB}}$ is indeed a Gröbner basis. Reducedness is obvious. Minimality is also obvious from Theorem 3.

## 6. Concluding remarks

In this paper we investigated the structure of degree two fibers of a decomposable model and clarified the structure of minimal Markov bases and minimal invariant Markov bases. We have also shown that decomposable models possess a Gröbner basis that is at the same time a minimal Markov basis.

For future research it is important to investigate structures of degree three fibers, degree four fibers, etc. In Takemura and Aoki (2004) we gave a characterization of minimal Markov bases. It shows that minimal Markov bases can be constructed "from below", that is, combining moves from fibers of degree $1, 2, 3, \ldots$. Although at the moment the construction cannot be implemented as an algorithm, it is important to study fibers of low degrees. We see that the study of degree two fibers in this paper led to some interesting results. As another example, in Aoki and Takemura (2009) we found some interesting degree three fibers in connection to experimental design with three-level factors.

As mentioned in the Introduction, the results obtained in this paper will provide insights to some practical models such as subtable sum models (Hara *et al.* (2009)), models for contingency tables with structural zeros (Aoki and Takemura (2005); Rapallo (2006)) and Rasch models (e.g., Chen and Small (2005); Zhu *et al.* (1998); Basturk (2008)) obtained by imposing some constraints on decomposable models. We will present results along this line in a forthcoming manuscript (Hara and Takemura (2009)).

It is of interest to study effects of structural zeros and observational zeros to Markov bases. In this respect in Hara *et al.* (2009) we have shown that a Markov basis for two-way contingency tables with structural zeros can be obtained as a subset of a Markov basis for subtable sum models, where the subtable sum happens to be an observational zero.

## Acknowledgements

The authors are grateful to Professor Takayuki Hibi of Osaka University and Professor Hidefumi Ohsugi of Rikkyo University for useful comments. The authors would also like to thank two anonymous referees for constructive comments and suggestions.